\documentclass[a4paper,11pt]{amsart}

\usepackage{amssymb}
\usepackage[pagebackref,
    ,pdfborder={0 0 0}
    ,urlcolor=black,hypertexnames=false]{hyperref}
\hypersetup{pdfauthor={Giuseppe Bargagnati, Alberto Casali, Francesco Milizia, Marco Moraschini},pdftitle=Minimal volume entropy of mapping tori}

\usepackage[all,cmtip]{xy}
\usepackage{appendix}
\usepackage{enumerate}
\usepackage{pgf}
\usepackage{tikz}
\usepackage{cleveref} 
\usetikzlibrary{decorations.pathmorphing,calc}
\tikzstyle{every picture}+=[remember picture]

\usepackage{tikz-cd}

\newtheorem{thm}{Theorem}[section]
\newtheorem{prop}[thm]{Proposition}
\newtheorem{lem}[thm]{Lemma}
\newtheorem{cor}[thm]{Corollary}

\newtheorem{quest}[thm]{Question}

\theoremstyle{remark}
\newtheorem{rem}[thm]{Remark}
\newtheorem{example}[thm]{Example}

\theoremstyle{definition}
\newtheorem{defi}[thm]{Definition}
\newtheorem{setup}[thm]{Setup}

\usepackage{amsfonts}
\newcommand{\Z}{\mathbb{Z}}

\newcommand{\R}{\mathbb{R}}
\newcommand{\N}{\mathbb{N}}

\newcommand{\connsum}{\mathbin{\#}}
\newcommand{\calP}{\mathcal{P}}
\newcommand{\calG}{\mathcal{G}}
\newcommand{\calF}{\mathcal{F}}

\DeclareMathOperator{\vol}{\textup{Vol}}
\DeclareMathOperator{\calU}{\mathcal{U}}

\DeclareMathOperator{\id}{id}

\DeclareMathOperator{\cat}{cat}
\DeclareMathOperator{\Am}{Am}
\DeclareMathOperator{\amcat}{\cat_{\Am}}

\DeclareMathOperator{\ent}{\textup{Ent}}
\DeclareMathOperator{\minent}{\textup{MinEnt}}
\DeclareMathOperator{\im}{im}

\DeclareMathOperator{\interior}{int}

\DeclareMathOperator{\Poly}{\textup{Poly}}
\DeclareMathOperator{\cut}{\backslash \! \backslash}

\DeclareMathOperator{\hol}{hol}
\DeclareMathOperator{\Hol}{Hol}

\DeclareMathOperator{\ess}{ess}

\makeatletter
\newcommand\norm{\bBigg@{0.8}}
\makeatother

\newcommand{\ifsv}[2][norm]{\csname #1l\endcsname\bracevert\!#2\!%
                            \csname #1r\endcsname\bracevert}
\newcommand{\ifsvp}[3][norm]{\csname #1l\endcsname\bracevert\!#2\!%
                            \csname #1r\endcsname\bracevert\!^{#3}}

\def\sv#1{\|#1\|}

\newcommand{\essn}[2][norm]{\csname #1l\endcsname\vert #2 \csname#1r\endcsname\vert_{1,\ess}}

\usepackage{color}
\usepackage{pdfcolmk}

\author[G.~Bargagnati]{Giuseppe Bargagnati}
\address{Dipartimento di Matematica\\ Universit\`{a} di Bologna\\ 40126~Bologna}
\email{giuseppe.bargagnati2@unibo.it}

\author[A.~Casali]{Alberto Casali}
\address{Dipartimento di Matematica \\ Universit\`{a} di Pisa \\ 56127~Pisa}
\email{alberto.casali@phd.unipi.it}

\author[F.~Milizia]{Francesco Milizia}
\address{Dipartimento di Matematica\\ Universit\`{a} di Bologna\\ 40126~Bologna}
\email{francesco.milizia@unibo.it}

\author[M.~Moraschini]{Marco Moraschini}
\address{Dipartimento di Matematica\\ Universit\`{a} di Bologna\\ 40126~Bologna}
\email{marco.moraschini2@unibo.it}


\title[Minimal volume entropy of mapping tori over 3-manifolds]{Minimal volume entropy of mapping tori over 3-manifolds}
\date{\today.\ \copyright{\ G.~Bargagnati, A.~Casali, F.~Milizia and M.~Moraschini}.}
\begin{document}

\begin{abstract}
  We prove that the minimal volume entropy of mapping tori over oriented closed smooth $3$-manifolds vanishes. Our approach uses a variation of the amenable category and a suitable version of the minimal volume entropy of a homology class introduced by Babenko and Sabourau.
\end{abstract}

\maketitle

\section{Introduction}

Given a closed Riemannian manifold~$(M, g)$, the \emph{volume entropy} of~$(M, g)$ measures the exponential growth rate of the volume of balls in the universal covering~$(\widetilde M, \tilde g)$ of $M$:
\[
\ent(M, g) := \lim_{R \to +\infty} \frac{1}{R} \cdot \log\big(\vol_{\tilde g}(B(x, R))\big),
\]
where $x$ is a point in the universal covering $\widetilde M$. A classical result of Manning~\cite{manning1979topological} ensures the existence of the limit and the independence on the chosen point $x$. This Riemannian invariant is related to the growth of the fundamental group by classical results of Svarc~\cite{vsvarc1955volume} and Milnor~\cite{milnor1968note} (in fact, it is positive if and only if the fundamental group of $M$ has exponential growth) as well as it is related to dynamics. Indeed, the volume entropy provides a lower bound for the topological entropy of the geodesic flow and the equality holds for non-positively curved manifolds~\cite{dinaburg1971connection, manning1979topological}. 

In order to obtain an invariant of the manifold $M$, starting from the volume entropy, one considers the infimum of $\ent(M, g)$ over suitable Riemannian metrics on $M$; an easy computation shows that $\ent(M, \lambda g) = \ent(M, g) \slash \sqrt{\lambda}$ for every $\lambda > 0$, justifying the following definition motivated by the work of Gromov~\cite{vbc} and then studied by many authors: The \emph{minimal volume entropy} of $M$ (known also with the name of \emph{absolute asymptotic volume}~\cite{babenko1993asymptotic}) is the real number
\begin{align*}
\minent(M) :&= \inf \{\ent(M, g) \vol(M, g)^{\frac{1}{\dim(M)}} \mid g \, \in\, \textup{Riem}(M)\} \\
&= \inf \{\ent(M, g) \mid {g \,  \in \, \textup{Riem}(M), \vol(M, g) = 1} \}.
\end{align*}
The minimal volume entropy is a homotopy invariant of the manifold $M$ and it only depends on the image of the fundamental class under the classifying map~\cite{brunnbauer2008homological}. Moreover, it is related to several other homotopy invariants, especially to \emph{simplicial volume}. Recall that the simplicial volume $\sv{M}$ of $M$ is defined as the infimum of all the $\ell^1$-norms of the real fundamental cycles of $M$, i.e., singular cycles representing the real fundamental class of $M$~\cite{vbc}. Gromov~\cite{vbc} established the following inequality:
\[
\minent(M)^{\dim(M)} \geq c_{\dim(M)} \cdot \sv{M}, 
\]
where $c_{\dim(M)}$ is a positive constant depending on the dimension of $M$ only. On the other hand, a natural classical question is whether it actually exists a biLipschitz inequality between the two invariants:

\begin{quest}\label{question:minent:vs:sv}
Let $M$ be an oriented closed connected smooth manifold; is it true that \[\sv{M} = 0 \Rightarrow \minent(M) = 0 \mbox{?}\]
\end{quest}

This question attracted a lot of attention but it is still widely open. The answer has long been known to be positive in dimension~$2$~\cite{vbc} and for the
class of manifolds whose fundamental group admits a representation with subexponential kernel into that of a negatively curved manifold~\cite{sambusetti1999minimal}. Moreover, it has been proved that an affirmative answer also holds in dimension~$3$~\cite{pieroni2019minimal} and $4$ provided that the $4$-manifold is geometrizable~\cite{suarez2009minimal}. On the other hand, Babenko and Sabourau have recently proved that the minimal volume entropy seminorm (a ``stabilization'' of the minimal volume entropy) is equivalent to the simplicial volume~\cite{BS2023seminorm}; this result can be understood as evidence towards a positive answer to Question~\ref{question:minent:vs:sv}.

One of the main difficulties in the computation of minimal volume entropy is that it is a quite crude invariant, and in particular it is not known whether it is multiplicative with respect to finite coverings~\cite{kedra2006crossed, kotschick2011entropies}.
In fact the situation is even worse, since it is not even known whether the vanishing of the minimal volume entropy of a finite covering ensures the vanishing of the minimal volume entropy of the base~\cite{kotschick2011entropies}.
It is possible that a counterexample to Question~\ref{question:minent:vs:sv} could be found by inspecting manifolds for which the simplicial volume is known to vanish because of the existence of finite coverings with zero simplicial volume.

In this paper we study the minimal volume entropy of all mapping tori over oriented closed smooth $3$-manifolds. Bucher and Neofytidis proved that the simplicial volume of such manifolds is always zero~\cite{bucherneofytidis} and their computation (especially for mapping tori over prime manifolds) strongly uses the multiplicativity of the simplicial volume under finite coverings. Nevertheless, we prove the following:

\begin{thm}[Theorem~\ref{thm:minent_reducible}]\label{thm:MT:prime:intro}
    The minimal volume entropy of all mapping tori of oriented closed smooth $3$-manifolds vanishes.
\end{thm}
We interpret the mapping tori in the smooth category, namely we take an oriented closed smooth $3$-manifold~$M$ and an orientation-preserving self-diffeomorphism $f \colon M \to M$; then, we denote by $M_f$ the mapping torus of $f$ (Remark~\ref{rem:MT:category:smooth:top}). This result shows that the answer to Question~\ref{question:minent:vs:sv} is affirmative also for this family of manifolds. 

One may wonder whether such mapping tori admit a geometry. This is already not true for some mapping tori over geometric prime $3$-manifolds~\cite[Theorem~9.8]{hillman2002four} (we refer the interested reader to the book by Hillman for a thorough discussion on geometric $4$-manifolds~\cite[Chapters~4, 8 and 9]{hillman2002four}). For this reason, our proof does not make explicit use of previous computations by Su\'arez-Serrato on the minimal volume entropy of oriented closed smooth $4$-manifolds admitting a geometry~\cite{suarez2009minimal}. 

In fact, Su\'arez-Serrato proved a stronger result on oriented closed smooth geometric $4$-manifolds~\cite[Theorem~A]{suarez2009minimal}: The vanishing of their minimal volume entropy is equivalent to the existence of an $F$-structure in the sense of Cheeger and Gromov~\cite{cheegergromov}. It is well-known that the existence of any $F$-structure on an oriented closed smooth manifold $M$ implies that its minimal entropy vanishes~\cite{paternainpetean};
hence, it is natural to ask whether every mapping torus of an oriented closed smooth $3$-manifold admits
such a structure. It is worth noting that in dimension $3$ mapping tori of surfaces generally
do not possess such a structure (specifically when $f$ is a pseudo-Anosov diffeomorphism of an oriented closed smooth surface of genus $g \geq 2$). At present, the authors are unable
to produce even a single example of a $4$-dimensional mapping torus without
an $F$-structure, so it might be possible that an analogue result of the one by Su\'arez-Serrato is also true in this setting. On the other hand, a natural source of possible counterexamples would be mapping tori over reducible $3$-manifolds, e.g., the mapping torus of a fully
irreducible, atoroidal diffeomorphism of $M = (S^2 \times S^1) \# (S^2 \times S^1)$. However, a formal proof of this non-existence result would require structural properties regarding the fundamental group of closed manifolds admitting an $F$-structure, which falls outside the scope of the present work.

\subsection{Strategy of the proof}

The proof is divided in two parts. We first prove the vanishing of the minimal volume entropy of mapping tori over oriented closed prime smooth $3$-manifolds and then we show that the result can be extended to all mapping tori over oriented closed smooth $3$-manifolds.
By contrast with the result on the simplicial volume by Bucher and Neofytidis, here also the computation of the minimal volume entropy of mapping tori over prime manifolds is a challenge.
Indeed, as stressed above we cannot rely, for our strategy of proof, on the existence of finite coverings with vanishing minimal volume entropy; instead, we have to prove directly that all the mapping tori have vanishing minimal volume entropy.
Our approach uses a variation of the \emph{amenable category} and a recent result by Babenko and Sabourau~\cite[Corollary~1.4]{BS2025fiber} (that, as already observed by Babenko and Saborau, can be explicitly stated in terms of covers~\cite[Corollary~5.9 and Example~5.10]{lmfibration}).
Recall that given a (possibly disconnected) open set~$U$ inside a manifold~$M$ we say that $U$ is \emph{amenable} (resp.\ \emph{of polynomial growth}) if the image of the group homomorphism induced by the inclusion
\[
\im \big(\pi_1(i) \colon \pi_1(U, x) \to \pi_1(M, x)\big)
\]
is amenable for all $x \in U$ (resp.\ such that all its finitely generated subgroups are of polynomial growth, i.e., finitely generated virtually nilpotent by a celebrated theorem of Gromov~\cite{gromov1981groups}). The \emph{amenable category of $M$}, denoted by $\amcat(M)$, is defined as the minimal integer $n \in \N$ such that there exists an open amenable cover of $M$ with cardinality $n$ (we stress that the open sets in $M$ are not assumed to be path-connected). 
The vanishing theorem of Gromov ensures that the simplicial volume of $M$ vanishes whenever $\amcat(M) \leq \dim(M)$~\cite[p.~40]{vbc}. An analogous result for the minimal volume entropy has been recently proved by Babenko and Sabourau~\cite[Corollary~1.4]{BS2025fiber}: The minimal volume entropy of $M$ vanishes whenever we can find an open cover of $M$ of cardinality $\leq \dim(M)$ consisting of \emph{sets of polynomial growth}. 
We denote by $\cat_{\Poly}(M)$ the minimum cardinality of such a cover. Hence, we prove the following result:

\begin{thm}[Theorem~\ref{thm:amcat:main}]\label{thm:amcat:intro}
    Let $M$ be an oriented closed prime smooth $3$-manifold and let $f \colon M \to M$ be an orientation-preserving self-\hspace{0pt}diffeomorphism. Then, 
    \[\amcat(M_f) \leq \cat_{\Poly}(M_f) \le 4 = \dim(M_f).\]
\end{thm}

This result is also of independent interest since it answers, for the class of manifolds given by mapping tori over oriented closed connected aspherical $3$-manifolds, the question whether the vanishing of the simplicial volume of an oriented closed connected aspherical manifold~$M$ implies $\amcat(M) \leq \dim (M)$~\cite[Question~1.3]{LMS} (a positive answer to this question for all closed connected aspherical manifolds would imply a positive answer to a celebrated question of Gromov about a relation between simplicial volume and Euler characteristic~\cite{LMS, loehmoraschiniraptis}).

Once we have obtained the vanishing of the minimal volume entropy for all mapping tori over oriented closed prime $3$-manifolds we have to extend this result to the case of mapping tori over \emph{all} oriented closed $3$-manifolds. Here, our approach follows the blueprint of the computation by Bucher and Neofytidis for the simplicial volume~\cite[Proof of Theorem~1.7]{bucherneofytidis}. However, also in this case the lack of multiplicativity with respect to finite coverings makes a straightforward adaptation not possible. Indeed, we have to work with 
a notion of minimal volume entropy of homology classes, introduced by Babenko and Sabourau~\cite{BS2023seminorm} (Section~\ref{Sec:minent:MT}).
This real-valued function on homology satisfies most of the required properties that Bucher and Neofytidis use about the $\ell^1$-norm~\cite[Proof of Theorem~1.7]{bucherneofytidis}, but it is still not multiplicative with respect to finite coverings. This forces us to diverge quickly from their proof and provide some new shortcuts that allow us to make the passage to finite coverings unnecessary. Such a strategy requires an accurate description of the maps and homotopies involved in the construction.

\subsection{Related results}

Given a group $\Gamma$ of type F (i.e., there exists a finite aspherical CW-complex of which it is the fundamental group), one can define the minimal volume entropy of $\Gamma$ as the infimum of the minimal volume entropies of finite aspherical simplicial complexes~$X$ of minimal dimension such that $\pi_1(X) \cong \Gamma$ \cite{bregman2021minimal}. The minimal volume entropy of a free group has been computed \cite{kapovich2007patterson,mcmullen2015entropy,Lim08}, 
and it is positive when the free group is non-Abelian, the infimum being realized by graphs whose vertices have degree $3$.
More recently, Bregman and Clay have provided lower estimates of the minimal volume entropy of free-by-cyclic groups~\cite[Theorem~1.1]{bregman2021minimal} and $2$-dimensional right-angled Artin groups~\cite[Theorem~1.2]{bregman2021minimal}. Also in this situation, a complete characterization of the case where the minimal entropy vanishes is possible. A subsequent natural question after the work by Bregman and Clay is to investigate the minimal volume entropy of extensions of $\mathbb{Z}$ by groups of higher cohomological dimension. 
Consider a short exact sequence~$1 \to N \to \Gamma \to \mathbb{Z} \to 1$, which corresponds to a semidirect product~$\Gamma = N \rtimes_{f_*} \mathbb{Z}$. Let $X$ be a cellular model for~$K(N, 1)$; then $f_*$ is induced by a self-homotopy equivalence $f \colon X \to X$ and we can thus identify the group~$\Gamma$ with the fundamental group of the mapping torus over~$X$ along the map~$f$. When $X$ is a closed manifold that satisfies the Borel conjecture, it follows that $f$ can be replaced up to homotopy by a self-homeomorphism of $X$; this applies to all oriented closed connected aspherical $3$-manifolds~\cite[Theorem~2.1.2]{aschenbrenner20123} as well as all oriented closed connected hyperbolic manifolds~\cite[Theorem~C.5.2]{benedetti1992lectures}. We say that the extension~$1 \to N \to \Gamma \to \mathbb{Z} \to 1$ is \emph{orientable} when the map~$f \colon X \to X$ is orientation-preserving. Theorem~\ref{thm:MT:prime:intro} together with Remark~\ref{rem:hyp:higher:dim} readily implies the following: 
%
%
%
%
%
%
%
%
%
\begin{cor}
    Let $N$ be the fundamental group of either an oriented closed connected aspherical $3$-manifold or an oriented closed connected hyperbolic manifold of dimension at least~$3$. Then every orientable group extension of~$\mathbb{Z}$ by~$N$ has vanishing minimal volume entropy.
%
%
\end{cor}

Another straightforward corollary of the previous results involves the stable integral simplicial volume and the integral foliated simplicial volume.
These invariants are suitable variations of classical simplicial volume either in terms of stabilization and integral cycles or in terms of dynamics~\cite{loehergodic, LMS}.
A classical question in this setting is whether oriented closed connected aspherical manifolds with residually finite fundamental group with vanishing simplicial volume also have zero integral foliated simplicial volume and stable integral simplicial volume~\cite[Question~6.2.2 and Question~6.3.3]{loehergodic}.
Since the fundamental group of mapping tori over oriented closed connected aspherical prime $3$-manifolds are residually finite~\cite[Theorem~II.44]{Steinhuber}, it is interesting to show that such manifolds satisfy the previous question. This was proved by Steinhuber~\cite{Steinhuber} using a careful adaptation of Bucher and Neofytidis techniques in the dynamical setting. Using Theorem~\ref{thm:amcat:intro} on the amenable category and the corresponding vanishing theorems for stable integral simplicial volume and integral foliated simplicial volume~\cite[Theorem~1.2]{LMS}, we can provide a new proof of this result (Remark~\ref{rem:hyp:higher:dim}):

\begin{cor}[{\cite[Theorem~V.I and Corollary~V.2]{Steinhuber}}]
    The stable integral simplicial volume and the integral foliated simplicial volume of all mapping tori over oriented closed connected aspherical prime $3$-manifolds vanish.

    The same conclusion is also true for mapping tori over oriented closed connected hyperbolic manifolds of dimension at least $3$.
\end{cor}

\subsection{Plan of the paper}

Section~\ref{Sec:amcat:glueings} is devoted to the study of the 
amenable category, and its analogous notion based on groups of polynomial growth, for manifolds obtained via glueing procedures as well as for manifolds admitting suitable foliations. This technical section is the heart of the computation of these invariants for mapping tori over oriented closed prime manifold, which is performed in Section~\ref{sec:amcat:mapping:tori:prime}, where we prove Theorem~\ref{thm:amcat:intro}. 
Finally, Section~\ref{Sec:minent:MT} contains the proof of Theorem~\ref{thm:MT:prime:intro} (Theorem~\ref{thm:minent_reducible}).
Appendix~\ref{sec:prelim:MT} contains the definition of mapping tori as well as many useful facts about their topology, homology and fundamental group, which are used throughout the paper.

\subsection{Acknowledgements}
We thank the referee for their careful proofreading and their report, which help to improve the exposition and the clarity of the paper.

This work has been funded by the European Union - NextGenerationEU under the National Recovery and Resilience Plan (PNRR) - Mission 4 Education and research - Component 2 From research to business - Investment 1.1 Notice Prin 2022 -  DD N.~104 del 2/2/2022, from title ``Geometry and topology of manifolds", proposal code 2022NMPLT8 - CUP J53D23003820001.

F.\ M.\ and M.\ M.\ were partially supported by the ERC ``Definable Algebraic Topology" DAT - Grant Agreement no.~101077154. 

All the authors also thank the INdAM-GNSAGA for supporting the conference ``\emph{Manifolds and groups in Bologna III}" in which parts of the paper have been discussed. 

This paper is part of the PhD project of Alberto Casali.

\section{$\calG$-category, glueings and foliations}\label{Sec:amcat:glueings}

In this section we prove some new results on the (generalized) amenable 
category of manifolds obtained from certain glueings and of manifolds 
supporting a smooth regular circle foliation. The combination of such results leads to Corollary~\ref{cor:final:glueing:foliations}, that is a fundamental step for the computation of the minimal volume entropy of mapping tori over oriented closed connected smooth prime $3$-manifolds in the next section. 

Let $\calG$ be an \emph{isq-class of groups}, i.e., a class of groups that is closed under isomorphisms, 
subgroups and quotients. The two main examples of $\calG$ that we consider in this paper are the 
class~$\Am$ of amenable groups and the class~$\Poly$ of those groups such that 
all their finitely generated subgroups have polynomial growth. With this notation, we can define $\calG$-sets and the $\calG$-category as follows: A possibly (possibly disconnected) open set~$U$ inside a manifold~$M$ is a $\calG$-\emph{set} if the image of the group homomorphism induced by the inclusion
\[
\im \big(\pi_1(i) \colon \pi_1(U, x) \to \pi_1(M, x)\big)
\]
lies in $\calG$ for all $x \in U$. Similarly, the $\calG$-\emph{category} of $M$ is the minimum cardinality of an open cover of $M$ consisting of (possibly disconnected)~$\calG$-sets. Of course, when $\calG = \Am$ the numerical invariant~$\cat_{\calG}(M)$ is the usual amenable category~$\amcat(M)$ of~$M$.

We will stress when the results only apply to the previous two classes of groups (namely, $\Am$ and $\Poly$) or to general isq-classes of groups.

\subsection{$\calG$-category and glueings}
Along this section $\calG$ denotes an isq-class of groups. We say that an open $\calG$-cover~$\calU$ (i.e., an open cover consisting of $\calG$-sets) of a compact connected manifold $M$ is \emph{efficient} if its cardinality is equal to $\cat_{\calG}(M)$.
A straightforward computation shows the following:

\begin{prop}[$\cat_{\calG}$ of manifolds glued along $\calG$-boundary components]\label{prop:glueings:amcat:standard}
Let $N$ be a compact (possibly disconnected) $n$-manifold and let $M$ be a compact (possibly disconnected) $n$-manifold obtained by glueing some pairwise homeomorphic boundary components of $N$ together. If all boundary components of $N$ that have been glued have fundamental group in the isq-class of groups~$\calG$ , then we have
\[
\cat_{\calG}(M) \leq \cat_{\calG}(N)+1.
\]
\end{prop}

\begin{rem}[glueing two connected manifolds together]
    We do not require $N$ to be connected, therefore the previous theorem can be applied to glueings between boundary components belonging to the same connected component, as well as different connected components. For instance a standard situation representing the latter case is when we glue together two compact connected manifolds along their homeomorphic (possibly disconnected) boundary; in this case $N$ is the disjoint union of the two manifolds and $M$ is the resulting manifolds after the glueing.
\end{rem}

\begin{proof}[Proof of Proposition~\ref{prop:glueings:amcat:standard}]
Set $K$ to be the union of all glueing loci in $M$ and let $H\cong K\times (-\varepsilon,\varepsilon)$ be a (possibly disconnected) bicollar neighbourhood of (the connected components of) $K$. By construction $H$ is a $\calG$-set of $M$. In order to extend this open $\calG$-set to an open $\calG$-cover of $M$, we have to give an explicit description of~$M$. Since the closure of $M\setminus H$ is homeomorphic to $N$, we can view $M$ as obtained by glueing the prescribed boundary components of~$N$ to the boundary of~$\overline{H}$. Set $S = N \cap \overline{H} \subset M$.
    Given an open efficient $\mathcal{G}$-cover~$\mathcal{U} =\{U_1, \cdots, U_{k}\}$ of $N$, we consider new open sets $\{U_1', \cdots, U_k'\}$ in $M$ given by \[U_i' = U_i \cup \Big((U_i\cap S)\times(-\varepsilon,-\varepsilon/2)\Big) \subset M = N \cup_S \overline{H}\]
    for every $i \in \{1, \cdots, k\}$.
    Then $\{U_1',\cdots, U_k'\}$ are $\calG$-sets in $M$ since each of them retracts to the corresponding $U_i$ (that is a $\calG$-set in $N \subset M$) and $\calG$ is closed under quotients. Thus we have provided an open $\calG$-cover of $M$ given by $\{H, U_1', \cdots, U_k'\}$. By taking the minimum among the cardinality of all the open $\calG$-cover of $M$ we get
    \[
    \cat_{\calG}(M) \leq \cat_{\calG}(N) + 1. \qedhere
    \]
\end{proof}

The previous inequality is sharp when $\calG$ is either $\Am$ or $\Poly$ as witnessed by the following example: 
\begin{example}[the $\cat_{\calG}$-inequality for glueings is sharp]
    Let $n \geq 2$ and let $M$ be an oriented compact connected $n$-manifold with non-empty boundary such that the interior of $M$ admits a complete hyperbolic metric of finite volume. Since $M$ is a compact manifold with non-empty boundary, it retracts by deformation onto a simplicial complex~$X$ of codimension one. In particular, this implies that~\cite[Remarks~2.8 and~2.9]{CLM} \[\cat_{\calG}(M) \leq \dim(X) + 1 \leq n.\]
    Let $D(M)$ denote the \emph{double} of $M$, that is the oriented closed connected $n$-manifold obtained by glueing two copies~$M_1, M_2$ of $M$ oriented in an opposite way along their boundary via the identity. Since we are glueing along $\pi_1$-injective virtually Abelian boundary components we have by the additivity of simplicial volume~\cite[Theorem~3]{bucher2014isometric}
    \[
    \sv{D(M)} = \sv{M_1} + \sv{M_2} = 2 \sv{M} > 0.
    \]
    By Gromov's vanishing theorem~\cite[p.~40]{vbc} we have that \[\cat_{\calG}(D(M)) > \dim(D(M)) = n = \cat_{\calG}(M).\] Hence, Proposition~\ref{prop:glueings:amcat:standard} implies 
    \[
    \cat_{\calG}(D(M)) = \cat_\calG(M)+1 = \cat_{\calG}(M \sqcup -M) + 1.\]
    
\end{example}

Even if in general we cannot hope to obtain a better estimate than the one in Proposition~\ref{prop:glueings:amcat:standard}, in many applications it would be very convenient to drop the $+1$ in the estimate above, and also to allow for more general glueings. This will be possible when the manifolds admit an efficient $\calG$-cover that is well behaved near the glueing site. To formalize this concept, we introduce the following definition:

\begin{defi}[$\mathcal{K}$-compatible $\calG$-cover]
    Let $M$ be a compact $n$-manifold and let $\mathcal{K} = \{K_1, \cdots, K_r\}$ be a collection of boundary components of $M$. We say that an open $\calG$-cover~$\calU$ of $M$ is $\mathcal{K}$-\emph{compatible} if the following holds: There exists an open set $U \in \calU$ such that 
    \begin{enumerate}
        \item It meets at least one boundary component (if any) of $M$ and no other open set $V \in \calU$ meets it (implying, in particular, that $U$ covers $\bigcup_{i = 1}^r K_i$);
        \item For every boundary component $K_i \in \mathcal{K}$ the path-connected component $U_i$ of $U$ intersecting $K_i$ is homeomorphic to a open collar neighbourhood~$(-\varepsilon, 0] \times K_i$.
    \end{enumerate}
\end{defi}
Henceforth, we suppose $(-\varepsilon, 0] \times K_i \subset M$ are fixed (parametrized) collar neighbourhoods.

\begin{rem}[cardinality of $\mathcal{K}$-compatible $\calG$-cover]\label{rem:K:compatible:cardinality:at:least:2}
    Item~$(2)$ implies that if $\mathcal{K} \neq \emptyset$ then an open $\mathcal{K}$-compatible $\mathcal{G}$-cover of a compact manifold must have cardinality at least $2$. 
\end{rem}

    Since each $K_i$ is a deformation retract of its open  collar neighbourhood, the existence of a $\mathcal{K}$-compatible $\calG$-cover implies that the boundary component $K_i$ is a (closed) $\calG$-set for every $i\in\{1,\cdots, r\}$. We can use this observation to give some equivalent definitions of $\mathcal{K}$-compatibility:
    \begin{prop}[equivalent definitions of $\mathcal{K}$-compatible $\calG$-cover]\label{prop:compatibility:equiv:defs}
        Let $M$ be a compact $n$-manifold, and $\mathcal{K}=\{K_1, \cdots, K_r\}$ be a non-empty collection of boundary components of $M$. 
        Then, the following conditions are equivalent:
        \begin{enumerate}
            \item There exists an open $\mathcal{K}$-compatible $\mathcal{G}$-cover $\mathcal{U}$;
            \item There exists an open $\mathcal{K}$-compatible $\mathcal{G}$-cover $\mathcal{U}$ such that for each path-connected component $P$ of every $V\in \calU$ the following holds: If the intersection~$P\cap K_i=\emptyset$, then $P\cap \big((-\varepsilon/2,0]\times K_i\big)=\emptyset$;
            \item For every $i\in\{1,\cdots, r\}$ the boundary component~$K_i$ is a $\calG$-set, and there exists an open $\mathcal{G}$-cover $\mathcal{U}$ of $M$ such that, for every $i$, every path-connected component $P$ of $V \in \mathcal{U}$ intersecting $K_i$ non-trivially is contained in the open collar neighbourhood $ (-\varepsilon, 0] \times K_i$.
        \end{enumerate}
    \end{prop}
    \begin{proof}
        Clearly  (\emph{2}) implies (\emph{1}) and (\emph{1}) implies (\emph{3}), therefore we only need to prove that (\emph{3}) implies (\emph{2}).
        
        Suppose that $\mathcal{U}$ satisfies the hypotheses in Item~(\emph{3}). Let us call $S = \cup_{i=1}^r K_i \subseteq \partial M$ and consider
        \[M' = M \cup_{S \cong S \times  \{0\}} \big(S \times [0, 1]\big),
        \]
        which is homeomorphic to $M$.
        Since $M$ is compact and $\mathcal{K}$ is non-empty, by Remark~\ref{rem:K:compatible:cardinality:at:least:2} we know that the  cardinality of $\calU$ is at least $2$. We construct an open $\calG$-cover~$\mathcal{V}$ of~$M'$ from $\calU = \{U_1, \cdots, U_k\}$, where $k\geq 2$, as follows: 
        \begin{align*}
        V_1 &= U_1 \cup \Big(\big(U_1 \cap S\big) \times [0, \varepsilon \slash 8)\Big) \cup \Big((\varepsilon \slash 4, 1] \times S)\Big), \\
        V_2 &= U_2 \cup \Big(\big(U_2 \cap S\big) \times [0, \varepsilon\slash 8)\Big) \cup \Big((0,\varepsilon \slash 2) \times S)\Big), \\
        V_j &= U_j \cup \Big(\big(U_j \cap S\big) \times [0, \varepsilon \slash 8)\Big)  \quad \quad \quad \quad  & \hspace{-50pt}\text{for every }j \in\{ 3, \cdots, r\}.
        \end{align*}
        All the path-connected components of the previous sets are $\calG$-sets in $M'$: Indeed, the only non-trivial computation appears in $V_2$, but its unique component intersecting $S$ is contained into $(-\varepsilon, \varepsilon \slash 2) \times S$ that is an open $\calG$-set by hypothesis. By construction, $\mathcal{V}$ is an open $\calG$-cover of $M'$ satisfying condition~(\emph{2}). By pulling back this open cover along the homeomorphism $M \cong M'$ we obtain the desired open $\calG$-cover of $M$ satisfying condition~(\emph{2}).
    \end{proof}

The notion of $\mathcal{K}$-compatible $\calG$-covers allows to prove the following glueing result:

\begin{prop}[$\cat_{\calG}$ of glueings in case of $\mathcal{K}$-compatible $\calG$-covers]\label{prop:glueing:general:K_i:new:version} Let $M$ be a compact $n$-dimensional manifold and let $K_1,K_2\subset\partial M$ be homeomorphic boundary components of $M$. Consider $f\colon K_1\to K_2$ a homeomorphism and denote with $M'$ the manifold obtained from $M$ by glueing $K_1$ to $K_2$ via~$f$.
Then, we have
\[
\cat_\calG(M')\leq \min \{|\calU| \ :\  \calU \text{ is } \{K_1,K_2\}\text{-compatible open }\mathcal{G}\text{-cover of } M\},
\]
with the convention that $\min \emptyset=+\infty.$

In particular, if $M$ admits a $\{K_1,K_2\}$-compatible efficient $\mathcal{G}$-cover, then
 \[
   \cat_\calG(M')\leq \cat_\calG(M).
    \] 
\end{prop}

\begin{proof}
    We prove the proposition in the case of $\{K_1,K_2\}$-compatible efficient $\mathcal{G}$-cover, the general proof being the same.
    Let $\calU$ be the given $\{K_1, K_2\}$-compatible efficient $\calG$-cover of $M$ and let $q \colon M \to M'$ be the quotient map corresponding to the glueing along the homeomorphism~$f \colon K_1 \to K_2$. By definition
    of $\{K_1, K_2\}$-compatibility, there exists a unique open set~$U \in \calU$ such that $U \cap \big(K_1 \cup K_2\big) \neq \emptyset$.
    We claim that the collection~$\mathcal{U}' = \{q(V)\}_{V \in \calU}$
    is an open $\calG$-cover of $M'$.
    Since for each $V \in \calU$ we have $q^{-1}(q(V)) = V$, it is immediate to check that $\mathcal{U}'$ is an open cover of $M'$.
    Moreover, since $\calU$ is $\{K_1, K_2\}$-compatible, we have that there is a unique path-connected component~$P'$ of a set in $\mathcal{U}'$ which intersects $q(K_1 \cup K_2)$; more precisely, this is a path-connected component of the set~$q(U)$. Moreover, because of $\{K_1, K_2\}$-compatibility of $\calU$ again, the path-connected component~$P'$ is homeomorphic to 
    $(-\varepsilon, \varepsilon) \times q(K_1 \cup K_2)$ (note that $q(K_1 \cup K_2) = q(K_1) = q(K_2) \cong K_1 \cong K_2$).
    So in order to show that $\mathcal{U}'$ is a $\calG$-cover we have to consider two distinct cases: If the path-connected component~$Q$ of some $q(V)$ is either $P'$ or not.

    We begin by assuming that the path-connected component~$Q' \neq P'$. In this case $Q'$
    is the homeomorphic image of the corresponding path-connected component~$Q$ of $V \in \calU$. Thus we have the following commutative diagram:
    \[
    \begin{tikzcd}
    \pi_1(Q) \ar[r] \ar[d] & \pi_1\left(M \setminus \big(K_1 \cup K_2\big)\right) \ar[r] \ar[d] & \pi_1(M) \ar[d]\\
    \pi_1(Q') \ar[r] & \pi_1\left(q\big(M \setminus \big(K_1 \cup K_2\big)\big)\right) \ar[r] & \pi_1(M'),
    \end{tikzcd}
    \]
    where all the vertical maps are induced by (restrictions of) $q$ and the horizontal maps by the inclusions. Since the two vertical arrows on the left are isomorphisms, the fact that $Q$ is a $\calG$-set in $M$ readily implies that also $Q' = q(Q)$ is a $\calG$-set in $M'$ (being $\calG$ closed under quotients).

    We are left to consider the path-connected component~$P'\subset q(U)$ that intersects $q(K_1 \cup K_2)$. By construction, the inclusion
    \[
    (-\varepsilon, -\varepsilon \slash 2) \times q(K_1 \cup K_2) \hookrightarrow P'
    \]
    is a homotopy equivalence, whence $P'$ is a $\calG$-set for $M'$ if and only if $(-\varepsilon, -\varepsilon \slash 2) \times q(K_1 \cup K_2)$ is so.
    This latter set is the homeomorphic image under $q$ of the set $(-\varepsilon, -\varepsilon \slash 2) \times K_1$ in $M$.
    Since this set does not intersect $K_1 \cup K_2$, as explained above, it is sufficient to show that this set is a $\calG$-set in $M$.
    Using the $K_1$-compatibility of the $\calG$-cover $\calU$ and the fact that the inclusion 
    \[
    (-\varepsilon, -\varepsilon \slash 2) \times K_1 \hookrightarrow (-\varepsilon, 0] \times K_1,
    \]
    is a homotopy equivalence,
    we readily see that $(-\varepsilon, -\varepsilon \slash 2) \times K_1$ is a $\calG$-set in $M$.
    Since the cardinality of $\mathcal{U}'$ is equal to the cardinality of $\mathcal{U}$, we have
    \[
    \cat_{\calG}(M') \leq \cat_{\calG}(M). \qedhere
    \]
\end{proof}

\begin{rem}[iterative construction]\label{rem:multiple:glueing:still:works}
  The proof of Proposition~\ref{prop:glueing:general:K_i:new:version} shows that given an open $\calG$-cover~$\calU$ of $M$, we do not modify the path-connected components of the sets $V \in \calU$ that do not intersect the boundary components $K_1$ and $K_2$ to construct the cover~$\calU'$ of~$M'$. This implies that if $\calU$ is a $\{K_1, K_2, \cdots, K_t\}$-compatible open $\calG$-cover of $M$, then the resulting open $\calG$-cover~$\calU'$ of $M'$ with $|\calU'| = |\calU|$ is still $\{K_3, \cdots, K_r\}$-compatible. Thus, we can iterate the construction and allow for multiple glueings in the statement of Proposition~\ref{prop:glueing:general:K_i:new:version} provided that the glueings are performed along pairwise-homeomorphic boundary components in $\{K_1, \cdots, K_r\}$.
\end{rem}

Let us change perspective now and assume that our manifolds are oriented. Let $M$ be an oriented compact connected $n$-manifold (with possibly empty boundary) and let $\mathcal{S} = \{S_1, \cdots, S_k\}$ be a collection of codimension-one closed connected disjoint submanifolds of $M \setminus \partial M$. Since each $S_i \in \mathcal{S}$ has codimension one, the tubular neighbourhood of $S_i$ in $M$ is simply given by an $I$-bundle. More precisely, if $S_i$ is orientable then the $I$-bundle is trivial (i.e., $S_i$ is two-sided), 
otherwise the $I$-bundle is twisted and its boundary is a connected double covering over $S_i$~\cite[p.~12]{hatcher2007notes} (i.e., $S_i$ is one-sided). Let us call $T(S_i)$ the tubular neigbourhood (possibly one-sided) of $S_i$ inside $M$. We consider the manifold
\[
M \cut \mathcal{S} := \overline{M \setminus (\cup_{i = 1}^k T(S_i))}.
\]
By construction, $M \cut \mathcal{S}$ is a disjoint union of oriented compact connected $n$-manifolds $M_1, \dots, M_t$ with non-empty boundaries (possibly with $t = 1$) and its interior is homeomorphic to $M \setminus (\cup_{i = 1}^k S_i)$. Moreover, we can reconstruct $M$ from $M \cut \mathcal{S}$ and the tubular neigbourhoods $T(S_i)$ by coeherently re-glueing them together. Thus from Proposition~\ref{prop:glueing:general:K_i:new:version} we obtain the following corollary:
\begin{cor}[$\calG$-estimate when cutting along submanifolds]\label{cor:cutting:amcat}
    Let $M$ be an oriented compact connected $n$-manifold. Let $\mathcal{S} = \{S_1, \cdots, S_k\}$ be a non-empty collection of codimension-one closed connected disjoint submanifolds of $M \setminus \partial M$, with fundamental groups $\pi_1(S_j) \in \calG$.
    Let $M_1, \dots, M_t$ be the connected components of the manifold~$M \cut \mathcal{S}$. Denote by $\mathcal{K}_i$ the set of the boundary components of $M_i$ that were not already boundary components of $\partial M$. Moreover,  suppose that there exists an integer $m \geq 2$ such that for every $i \in \{1, \cdots, t\}$  there exists and open  $\calG$-cover~$\mathcal{U}_i$ of $M_i$ that is $\mathcal{K}_i$-compatible and has cardinality $|\mathcal{U}_i|\leq m$. Then we have 
    \[
    \cat_\calG(M)\leq m.
    \]
    
    In particular, if for every $i \in \{1, \cdots, t\}$  there exists and open efficient $\calG$-cover~$\mathcal{U}_i$ of $M_i$ which is $\mathcal{K}_i$-compatible, then
    \[
    \cat_{\calG}(M) \leq \max \{\cat_{\calG}(M_1), \cdots, \cat_{\calG}(M_t)\}.
    \]
\end{cor}
\begin{proof}
    It is sufficient to show that under the hypotheses the disconnected compact manifold
    \[
    M' = \bigsqcup_{i = 1}^t M_i \sqcup \bigsqcup_{i = 1}^k T(S_i),
    \]
    has an open efficient $\{\mathcal{K}_1, \cdots, \mathcal{K}_t,\partial T(S_1),\cdots, \partial T(S_k)\}$-compatible $\calG$-cover $V$ of cardinality at most $m$  and to apply Proposition~\ref{prop:glueing:general:K_i:new:version} and Remark~\ref{rem:multiple:glueing:still:works}.

    For each tubular neighbourhood~$T(S_i)$, we construct a boundary-\hspace{0pt}compatible open $\calG$-cover of cardinality two.
    As mentioned above we have two different types of tubular neighbourhoods according to whether $S_i$ is one-sided or two-sided.
    If $S_i$ is two-sided we have that 
    \[
    T(S_i) \cong [-\varepsilon, \varepsilon] \times S_i.
    \]
    Hence, using the fact that $\pi_1(S_i) \in \calG$ it is immediate to check that 
    \[
    \left\{\left(\left[-\varepsilon, -\frac{\varepsilon}{4}\right) \cup \left(\frac{\varepsilon}{4}, \varepsilon\right]\right) \times S_i, \left(-\frac{\varepsilon}{2}, \frac{\varepsilon}{2}\right)\times S_i\right\}
    \]
    is an open boundary-compatible $\calG$-cover of $T(S_i)$.
    
    On the other hand, if $S_i$ is one-sided,  we consider the corresponding double covering~$\partial T(S_i) \to S_i$ with $\mathbb{Z} \slash 2 \Z = \{[0], [1]\}$ as deck transformations. Then, we have
    \[
    T(S_i) \cong \frac{[-1, 1] \times \partial T(S_i)}{(a, x) \sim (-a, [1]\cdot x)}.
    \]
    Hence, an open boundary-compatible $\calG$-cover of $T(S_i)$ is given by the following two connected open sets:
    \[
    \left\{\frac{[-1, -\frac{1}{4}) \cup (\frac{1}{4}, 1] \times \partial T(S_i)}{(a, x) \sim (-a, [1]\cdot x)},\frac{(-\frac{1}{2}, \frac{1}{2}) \times \partial T(S_i)}{(a, x) \sim (-a, [1]\cdot x)} \right\},
    \]
    where the first set is homotopy equivalent to $\partial T(S_i)$ and the second set is homotopy equivalent to $S_i$. Since $\pi_1(S_i) \in \calG$ by hypothesis and $\pi_1(\partial T(S_i))$ is an index-two subgroup of $\pi_1(S_i)$, we have that both sets are $\calG$-sets in $T(S_i)$ because $\calG$ is closed under subgroups.

    We have thus proved that all the connected components of $M'$ admit an open $\calG$-cover compatible with the desired set of boundary components (that is of cardinality at most $m \geq 2$ for each~$M_i$ and of cardinality two for each tubular neighbourhood~$T(S_i)$). We thus obtain by Proposition~\ref{prop:glueing:general:K_i:new:version} that
    \[
    \cat_{\calG}(M) \leq m.
    \]
    If every $M_i$ admits an efficient $\mathcal{K}_i$-compatible $\mathcal{G}$-cover, we can take
    \[
    m=\max \{\cat_{\calG}(M_1), \cdots, \cat_{\calG}(M_t)\}.\qedhere
    \]
\end{proof}

\begin{rem}[dropping the assumption of $\pi_1(S_i) \in \calG$]
    In fact, using an argument similar to the one in the proof of Proposition~\ref{prop:glueing:general:K_i:new:version}, one can drop the assumption of $\pi_1(S_i) \in \calG$ from the previous result. However, since we need this variant only, we preferred to restrict to the special situation of Corollary~\ref{cor:cutting:amcat}.
\end{rem}

\subsection{$\calG$-category and foliations}
$\calG$-covers with small multiplicity often arise from geometric situations, for instance in presence of certain foliations (Proposition~\ref{prop:amcat:circle:foliation:with:boundary}). 
For the computation of the minimal volume entropy of mapping tori over oriented closed connected smooth prime $3$-manifolds it will be crucial to combine the glueing techniques developed in the previous sections with the existence of amenable (or, better, of polynomial growth) covers with small cardinality associated to smooth regular circle foliations (Corollary~\ref{cor:final:glueing:foliations}). Recall that a \emph{smooth regular foliation} of a smooth manifold~$M$ is a partition~$\{L_p \, | \, p \in M\}$ of $M$ in \emph{leaves} such that
\begin{itemize}
    \item Each leaf~$L_p$ is an injectively immersed smooth $q$-submanifold of $M$ (called \emph{leaf through the point $p$});
    \item The tangent spaces of the leaves, taken together, form a smooth (integrable) subbundle of the tangent space~$TM$ of $M$.
\end{itemize}
Given a smooth regular foliation~$\mathcal{F}$ of a compact smooth manifold~$M$ with non-empty boundary~$\partial M$, we say that $\partial M$ is \emph{vertical} if it is union of leaves of $\mathcal{F}$ (this terminology is not standard). For convenience when $\partial M = \emptyset$ we say that $\partial M$ is vertical. A recurring example in this paper is the following:

\begin{example}[Foliation by suspension~{\cite[p.~16]{moerdijk-book}}]\label{ex:foliation:suspension}
    Let $M$ be an oriented compact connected smooth manifold with (possibly empty) boundary and let $f \colon (M, \partial M) \to (M, \partial M)$ be an orientation-preserving self-diffeomorphism of pairs. Then, we can identify the mapping torus~$M_f$ with the following quotient space:
    \[
    M_f \cong \frac{M \times \R}{(x, t) \sim (f^k(x), t - k)},
    \]
    for all $k \in \mathbb{Z}$ and all $(x, t) \in M \times \R$. It is immediate to check that $M \times \R$ has a foliation~$\mathcal{F}$ in lines where the leaves are simply $\{x\} \times \R$.
    Moreover, the properly discontinuous action of $\Z$ sends leaves to leaves, so we obtain a new foliation~$\mathcal{F}_f$ on the mapping torus~$M_f$ given by quotients of lines. By construction the boundary~$\partial M_f$ is vertical with respect to $\mathcal{F}_f$. This foliation is called a \emph{foliation by suspension}.

    Moreover, if $M$ is already foliated by a foliation~$\mathcal{F}$ and $f$ maps leaves to leaves, we obtain that $M \times \R$ is foliated by a foliation whose leaves are $\{F \times \R\}_{F \in \mathcal{F}}$. One can check that the $\Z$-action sends leaves to leaves again and so $M_f$ is foliated by quotients of $F \times \R$. If $\partial M$ is vertical with respect to the initial foliation~$\mathcal{F}$, then also $\partial M_f$ is vertical with respect to the new induced foliation by suspension.
    
\end{example}
A useful tool used to describe the local behaviour of a foliation is the notion of \emph{holonomy}. We will briefly recall its definition, for more details we refer the reader to the book by Moerdijk and Mr\v{c}un~\cite[Section~2.1]{moerdijk-book}. Intuitively, the holonomy is a representation of $\pi_1(L_p)$ that describes what happens if we ``follow'' the foliation on transverse submanifolds along a curve. 
Formally, for every $x\in M$ let $U_x$ be a foliated neighbourhood, i.e., an open neighbourhood of $x$ with a diffeomorphism $\rho\colon U_x \to \R^k \times \R^{n-k}$ such that $\rho(x)=(0,0)$ and that $\mathcal{F}$ restricted to $U_x$ is of the form $\{\R^k\times \{y\} \}_{y\in \R^{n-k}}$.
Consider two points $p,q\in L_x\cap U_x$ lying in the same path-connected component (namely, in the same \emph{plaque}) and fix $S$ and $T$ two submanifolds of $M$ transverse to $\calF$ in $p$ and $q$, respectively.
This means that, up to restricting $U_x$, the images $\rho(S\cap U_x)$ and $\rho(T\cap U_x)$ are of the form $\{(s(z),z)\}_{z\in \R^{n-k}}$ and $\{(t(z),z)\}_{z\in \R^{n-k}}$ for some smooth maps~$s,t\colon \R^{n-k}\to \R^k$, respectively.
Now the map  $(s(z),z)\mapsto (t(z),z)$ induces a diffeomorphism $f\colon S\cap U_x \to T\cap U_x$ sending $p$ to $q$.
Then we can consider the \emph{germ} of the map $f$ at $p$
\[ 
\hol^{S,T}:= \textup{germ}_p (f) \colon (S,p)\to (T,q).
\]
Recall that the \emph{germ} of a map $f$ at $p$ is the equivalence class of maps that coincide with $f$ on a neighbourhood of $p$.

Consider now a leaf $L$ and a closed path $\alpha\colon [0,1]\to L$ based in $p\in L$. We can find a sequence of times $t_0=0< t_1<\dots<t_m=1$ such that $\alpha(t_{i+1})$ lies in $U_{\alpha(t_i)}$. We can then consider a sequence $T_0, T_1, \dots, T_m=T_0$ of submanifolds such that each $T_i$ is transverse to $\calF$ in $\alpha(t_i)$. We can then define the holonomy of $\alpha$ as
\[ \hol^{T_0}(\alpha)= \hol^{(T_{m-1},T_m)}\circ \dots \circ \hol^{(T_{1},T_0)}\colon (T_0,p)\to (T_0,p). \]

One can prove that  $\hol^{T_0}(\alpha)$ depends only on the homotopy class of $\alpha$ in $L$ and the composition of loops is sent to the composition of the relative germs~\cite[p.~22]{moerdijk-book}. Hence, this procedure defines a group homomorphism
\[\hol^{T_0}\colon \pi_1(L,p) \to \textup{Diff}_p(T_0),\]
where $\textup{Diff}_p(T_0)$ is the group of germs of local diffeomorphisms of $T_0$ at $p$.
 Since $\textup{Diff}_p(T_0)$ is isomorphic to $\textup{Diff}_0(\R^{n-k})$, we can define a group homomorphism
\[\hol\colon \pi_1(L,p)\to \textup{Diff}_0(\R^{n-k})\]
that does not depend, up to conjugation in $\textup{Diff}_0(\R^{n-k})$, on the choice of the transverse submanifold $T_0$ and of the basepoint $p\in L$ and it is called the \emph{holonomy homomorphism} of $L$. The image of the holonomy homomorphism is called the \emph{holonomy} of the leaf $L$, and it is denoted by $\Hol(L)$~\cite[p.~23]{moerdijk-book}.
A leaf is said to be \emph{ordinary} if it has trivial holonomy, otherwise it is said to be \emph{exceptional}. The following classical result by Reeb ensures that in case of finite holonomy the neighbourhood around the leaf is foliated and well-behaved. 

\begin{thm}[{Local Reeb stability theorem~\cite[Theorem~2.9 and Remark below it]{moerdijk-book},~\cite[Theorem~3]{camachoneto}}] \label{thm:reeb:stability}
    Let $M$ be an oriented closed connected smooth $n$-manifold and $\mathcal{F}$ be a smooth regular foliation in codimension-$q$ submanifolds of $M$.
    Let $L \in \mathcal{F}$ be a compact leaf with finite holonomy group~$H$.
    Then there exists a foliated (i.e., union of leaves) arbitrarily small open set~$V_L$ containing $L$ and which is diffeomorphic to
    \[
    \overline{L} \times_H \R^{q} := \overline{L} \times \R^q \slash (l h, x) \sim (l, h x),
    \]
    where $\overline L \to L$ is the holonomy covering, i.e., the covering corresponding to the kernel of the holonomy homomorphism, and $H$ acts on $\overline{L}$ by deck transformations and on $\R^{q}$ linearly.
\end{thm}

\begin{rem}[local description in the case of compact manifolds with boundary]\label{rem:reeb:boundary}
    Theorem~\ref{thm:reeb:stability} can be extended to the case of an oriented compact connected smooth manifold~$M$ with non-empty vertical boundary by taking the double~$D(M)$ of~$M$.
    More precisely, every compact leaf~$L\subset\partial M$  with finite holonomy admits a foliated neighbourhood~$V_L$ as above in $D(M) = M \cup_{\partial M} (-M)$. The intersection $V_L \cap M$ provides the foliated neighbourhood of $L$ for the compact manifold with boundary $M$.
\end{rem}

We begin by extending the computation of the amenable category of closed smooth manifolds foliated by circles with finite holonomy~\cite[Proposition~7.3]{LMS} to the case of the $\calG$-category of compact smooth manifolds, and allowing for foliations of higher dimension:

\begin{prop}[$\cat_{\calG}$-estimate for compact manifolds foliated by certain compact leaves]\label{prop:amcat:circle:foliation:with:boundary}
    Let $\calG$ be an isq-class of groups. Let $M$ be an oriented compact connected smooth $n$-manifold~$M$ with (possibly empty) boundary~$\partial M$. Suppose that $M$ admits a smooth regular  foliation~$\mathcal{F}$ of codimension $q\leq n-1$. Suppose that every $L\in \calF$ is compact, has finite holonomy and $\pi_1(L)\in\calG$ and assume  that $\partial M$ is vertical with respect to $\mathcal{F}$. If $\partial M$ is a $\calG$-set, then there exists an open $\calG$-cover~$\calU$ of $M$ that is $\partial M$-compatible and has cardinality $q+1 \leq \dim(M)$. 
\end{prop}
\begin{proof}
    We follow the proof for circle foliations in the closed case as blueprint~\cite[Proposition~7.3]{LMS} and we start by treating the case where $\partial M\neq \emptyset$. Let $D(M)$ be the double of $M$. We are going to construct a nice open $\calG$-cover of $D(M)$ of cardinality $n$ and then we modify it in order to obtain the desired $\partial M$-compatible open $\calG$-cover of $M$.

    As discussed in Remark~\ref{rem:reeb:boundary} since the boundary~$\partial M$ is vertical, $D(M)$ admits a foliation $D(\mathcal{F)}$ that by construction satisfies the hypotheses. Denote by $\tau\colon D(M)\to D(M)$ the natural reflection along $\partial M$. Then by construction $\tau$ preserves the foliation $D(\mathcal{F})$. Let us denote by~$X$ the leaf space of $D(\mathcal{F})$ endowed with the quotient topology, i.e., $X=D(M)/D(\mathcal{\calF})$. Let $\pi\colon D(M)\to X$ be the projection and for every
    leaf~$L\in D(\mathcal{F})$, we denote by~$V_L$ the foliated open neighbourhood described in Theorem~\ref{thm:reeb:stability}. Since $D(\calF)$ is symmetric we can choose each $V_L$ such that $\tau(V_L)=V_{\tau(L)}$. Then the collection of $\pi(V_L)\cong \R^{q}/H$ provides a smooth orbifold atlas for $X$ of dimension $\dim(X)=q$~\cite[p.~40]{moerdijk-book}. In particular, $X$ is triangulable~\cite[Proposition~1.2.1]{moerdijk+pronk}. Now fix a bicollar neighbourhood $(-\varepsilon, \varepsilon)\times \partial M$ of $\partial M$ in $D(M)$ that is $\tau$-invariant. Up to choosing a smaller foliated open neighbourhood~$V_L$ of $L$ (Theorem~\ref{thm:reeb:stability}), we can assume that for every leaf~$L\subset \partial M$ we have that $V_L\subset(-\varepsilon, \varepsilon)\times \partial M$, and for every leaf~$L\in D(M)\setminus \partial M$ we have that $V_L\cap \partial M=\emptyset$. Since $D(M)$ is compact, $X$ is also compact and we can find a finite collection of leaves~$\{L_1,\cdots, L_k\}$ such that $\cup_{i=1}^k \pi(V_{L_i})=X$. Thus, up to taking a finite amount of barycentric subdivisions, we can assume that the \emph{open star}~$U_v$ of each vertex of~$X$ is contained in some open set~$\pi(V_{L_i})$. Moreover, if we assume that we performed at least one barycentric subdivision, every vertex of the subdivision corresponds to a simplex of the previous triangulation, and the open stars of two vertices corresponding to simplices of the same dimension (in the previous triangulation) are disjoint. 
    Therefore, if we denote by $U_i$ the union of all the open stars~$U_v$ where $v$ is a vertex in the last triangulation corresponding to an $i$-dimensional simplex of the penultimate triangulation, we have that $\{U_i\}_{i=0}^{q}$ is an open cover of $X$ of cardinality $n$ (because $\dim(X) = q$). 
    We claim that $\{\pi^{-1}(U_i)\}_{i=0}^{q}$ is an open $\mathcal{G}$-cover of $D(M)$ of cardinality $q+1$. Since each each~$\pi^{-1}(U_v)$ is a subset of some $V_L$ and $\calG$ is closed under subgroups and quotients, it is sufficient to notice that $V_L$ is a $\calG$-set for every leaf~$L$.
    This readily follows from the fact that $V_L\cong \overline{L}\times_H \R^{q} \subset D(M)$ has fundamental group $\pi_1(V_L)\cong \pi_1(L) \in \calG$ because it is a fiber bundle over $L$ with fiber $\R^q$~\cite[p.~17]{moerdijk-book}. This proves the claim. We are now going to produce the desired open $\partial M$-compatible open $\calG$-cover of $M$ of cardinality $q+1$ from this cover.
    


    Consider the open cover~$\mathcal{U}=\{\pi^{-1}(U_i) \cap M\}_{i=0}^{q}$. To show that $\mathcal{U}$ is an open $\mathcal{G}$-cover of $M$ it is sufficient to show that for each leaf~$L\in \calF$ the open set~$V_L\cap M$ is a $\calG$-set in $M$ (because $\calG$ is closed under quotients and subgroups). We have two cases:
    If $L\subset M\setminus\partial M$, then by construction $V_L\cap M=V_L$, hence $V_L \cap M$ is a $\calG$-set in M. Otherwise, if $L\subset\partial M$, we can consider the retraction $r\colon D(M) = M \cup -M \to M$ defined by
   \[
   r(m)=\begin{cases}
       m &\text{if } m\in M\\
       \tau(m) &\text{if } m\in -M.
   \end{cases}
   \]
    Since by construction we chose $V_L$ such that $\tau(V_L)=V_{\tau(L)}=V_L$, the restriction of $r$ to $V_L$ is a retraction of $V_L$ onto $V_L\cap M$.
    Thus, $\pi_1(V_L\cap M)$ is a quotient of $\pi_1(V_L)\in \calG$ and so, since $\calG$ is closed under quotients, we have that $\pi_1(V_L\cap M)\in \calG$. 

    We are left to show that $\calU$ is $\partial M$-compatible. Because of Proposition~\ref{prop:compatibility:equiv:defs}, it is sufficient to prove that~$\calU$ has the following property: For every path-connected component $P$ of $\pi^{-1}(U_i) \cap M \in \calU$ such that $P \cap \partial M \neq \emptyset$ we have that $P$ is contained in the open collar neighbourhood~$(-\varepsilon, 0] \times \partial M$. Suppose that $\pi^{-1}(U_v) \cap \partial M \neq \emptyset$ for some vertex $v$ in the refined triangulation of $X$. Then, by construction also the corresponding open set~$V_{L_i}$ (i.e., the set such that $\pi(V_{L_i}) \supset U_v$) intersects the boundary~$\partial M$ non-trivially. Thus, because of our choice of the foliated open sets~$V_L$ we have
    \[
    \pi^{-1}(U_v) \subset V_{L_i} \subset (-\varepsilon, \varepsilon) \times \partial M \subset D(M).
    \]
    This shows that the intersection $\pi^{-1}(U_v) \cap M$ is entirely contained in the collar~$(-\varepsilon, 0] \times \partial M$ and concludes the proof in the case when $\partial M\neq \emptyset$.

    The proof in the closed case is exactly the same, but instead of using the double $D(M)$ one can work directly with $M$ and the projection map $\pi\colon M\to M/\calF$. In this case, the pull-back of open cover $\{\pi^{-1}(U_i)\}_{i=0}^q$ is already the desired $\calG$-cover of cardinality $q+1$.
%
 %
\end{proof}

We then obtain the following corollary:

\begin{cor}[$\cat_{\calG}$-estimate for glueings of manifolds foliated by circles]\label{cor:final:glueing:foliations}
    Let $\calG$ be either $\Am$ or $\Poly$.
    Let $M$ be an oriented closed connected smooth $n$-manifold.
    Let $\mathcal{S} = \{S_1, \cdots, S_k\}$ be a non-empty collection of codimension-one closed connected disjoint smooth submanifolds of $M \setminus \partial M$, whose fundamental groups $\pi_1(S_j)$ are in $\calG$.
    Let $M_1,\dots,M_t$ be the connected components of the smooth manifold~$M \cut \mathcal{S}$. Suppose that for every $i \in \{1, \cdots, t\}$ the manifold~$M_i$ admits a smooth regular circle foliation~$\calF$ with finite holonomy such that the boundary $\partial M_i$ is vertical with respect to $\calF$. Then, we have:
    \[
    \cat_{\calG}(M) \leq n = \dim(M).
    \]
\end{cor}
\begin{proof}
   It is sufficient to apply Proposition~\ref{prop:amcat:circle:foliation:with:boundary} and Corollary~\ref{cor:cutting:amcat}.
\end{proof}

\section{$\calG$-category of mapping tori over prime smooth $3$-manifolds}\label{sec:amcat:mapping:tori:prime}

In this section, we consider the isq-class of groups~$\calG$ to be either $\Am$ or $\Poly$ and we compute the $\calG$-category of mapping tori over closed prime $3$-manifolds:
\begin{thm}[$\cat_{\calG}$ of mapping tori over prime $3$-manifolds]\label{thm:amcat:main}
   Let $M$ be an oriented closed prime smooth $3$-manifold and let $f \colon M \to M$ be an orientation-preserving self-diffeomorphism. Then, 
    \[
    \amcat(M_f) \leq \cat_{\Poly}(M_f) \leq \dim(M_f) = 4.
    \]
%
\end{thm}
Using a recent result by Babenko and Sabourau~\cite[Corollary~1.4]{BS2025fiber} this theorem readily implies the vanishing of the minimal volume entropy:
\begin{cor}[$\minent$ of mapping tori over prime $3$-manifolds]\label{cor:minent:MT:prime}
   Let $M$ be an oriented closed prime smooth $3$-manifold and let $f \colon M \to M$ be an orientation-preserving self-diffeomorphism. Then, 
    \[
    \minent(M_f) = 0.
    \]
\end{cor}
\begin{proof}
    Since by Theorem~\ref{thm:amcat:main} $\cat_{\Poly}(M_f) \leq \dim(M_f)$, we have that there exists an open cover of $M_f$ of cardinality less than or equal to $\dim(M_f)$ consisting of sets of polynomial growth. This is equivalent to the fact that there exists an iterated barycentric subdivision of $M_f$ that satisfies the fibre collapsing assumption with polynomial growth~\cite[Example~5.10]{lmfibration}. Since the minimal volume entropy is a homotopy invariant (and so, in particular, invariant under homeomorphisms), this implies that the minimal volume entropy of $M_f$ vanishes~\cite[Corollary~1.4]{BS2025fiber}.
\end{proof}

The amenable category of closed $3$-manifolds has been successfully computed by G\'omez-Larra{\~n}aga, Gonz\'alez-Acu\~na and Heil~\cite{GGH} and so it is possible to estimate the $\calG$-category of the mapping tori via the following formula:
\begin{thm}[$\cat_{\calG}$-estimate for mapping tori~{\cite[Theorem~1.1]{lmfibration}}]\label{thm:mapping:tori:classical:estimate:LM}
    Let $\calG$ be an isq-class of groups.
    Let $M$ be an oriented closed smooth manifold and let $f \colon M \to M$ be an orientation-preserving self-diffeomorphism of $M$. Then, we have
    \[
    \cat_{\calG}(M_f) \leq 2 \cdot \cat_{\calG}(M).
    \]
\end{thm}
However, already for estimating the amenable category of a $4$-dimensional mapping torus the previous formula is very loose as witnessed, e.g., for mapping tori over hyperbolic $3$-manifolds (in that case the estimate $\leq 8$ is even worse than the usual estimate $\leq \dim(M_f) + 1$). 

The goal of this section is to prove that one can in fact obtain the estimate
\[
\cat_{\calG}(M_f) \leq \dim(M_f) = 4
\]
for all oriented closed prime smooth $3$-manifolds~$M$ and all orientation-preserving self-diffeomorphisms~$f \colon M \to M$. This will give a proof of Theorem~\ref{thm:amcat:intro}. Before discussing all the different situations, recall that a group~$\Gamma$ is $P$\emph{-by-cyclic} if it is an extension of a group $\Lambda$ with property $P$ by the \emph{infinite} cyclic group~$\Z$. Thus, given a closed manifold~$M$ whose fundamental group has a property~$P$, then the fundamental group of every mapping torus over $M$ with respect to homeomorphism is $P$-by-cyclic (Remark~\ref{rem:mapping:torus:group:semidirect}).

In the subsequent discussion, we will use extensively the Geometrization Theorem for closed 3-manifolds~\cite{Per02, Per03bis, Per03, CaoZhu, MorganTian, KleLott, francesi}.

\subsection{Mapping tori over non-aspherical prime smooth $3$-manifolds}
Let us consider oriented closed connected non-aspherical prime $3$-manifolds. In this situation the manifold~$M$ is either \emph{elliptic} (i.e., covered by $S^3$) or diffeomorphic to $S^2 \times S^1$. In both cases the computation of $\cat_{\calG}$ is immediate: 
\begin{prop}[$\cat_{\calG}$-estimate for mapping tori over non-aspherical prime $3$-manifolds]\label{prop:MT:non:aspherical}
    Let $\calG$ be either $\Am$ or $\Poly$.
    Let $M$ be an oriented closed non-aspherical prime smooth $3$-manifold and let $f \colon M \to M$ be an orientation-preserving self-\hspace{0pt}diffeomorphism. Then, we have
    \[
    \cat_{\calG}(M_f)=1.
    \]
\end{prop}
\begin{proof}
   We have two cases. When $M$ is elliptic, the fundamental group~$\pi_1(M)$ is finite and so the fundamental group of the mapping torus~$\pi_1(M_f)$ is finitely generated finite-by-cyclic, thus finitely generated virtually cyclic~\cite[Lemma~2.2]{pettet1995finitely}. This shows that $\pi_1(M_f) \in \Poly \subset \Am$ and so $\cat_{\calG}(M_f) = 1$. On the other hand, if $M = S^2 \times S^1$ the fundamental group is infinite cyclic, and so $\pi_1(M_f)$ is cyclic-by-cyclic.
   Since $|\textup{Aut}({\Z})|=2$,  $\pi_1(M_f)$ is either isomorphic to $\Z\times \Z$ or to the non-Abelian semidirect product $\Z\rtimes \Z$ associated to the alternating homomorphism $\Z\to \textup{Aut} (\Z)$.
   In this second case, we have an index $2$ subgroup $\Z \rtimes 2\Z\cong \Z\times\Z$. This shows that $\pi_1(M_f)$ is finitely generated virtually Abelian and so $\pi_1(M_f) \in \Poly \subset \Am$. This proves that $\cat_{\calG}(M_f) = 1$.
\end{proof}

\subsection{Mapping tori over amenable aspherical prime smooth $3$-\hspace{0pt}manifolds}

Let $M$ be an oriented closed connected aspherical prime $3$-\hspace{0pt}manifold that has infinite a\-me\-na\-ble fundamental group. This means that $M$ admits one of the following three geometries: $\mathbb{E}^3, \textup{Nil}$ and $\textup{Sol}$. 
\begin{prop}[$\cat_{\calG}$-estimate for mapping tori over amenable aspherical prime $3$-manifolds]\label{prop:MT:infinite:amenable}
   Let $M$ be an oriented closed aspherical prime smooth $3$-manifold with infinite amenable group and let $f \colon M \to M$ be an orientation-preserving self-diffeomorphism. Then, we have 
   \[
   1 = \amcat(M_f) \leq \cat_{\Poly}(M_f) \leq 4.
   \]
\end{prop}
\begin{proof}
    By construction the fundamental group of $M_f$ is amenable-by-cyclic, thus amenable~\cite[Proposition~3.4]{frigerio2017bounded}. This shows that $\amcat(M) = 1$. On the other hand, it is not true in general that amenable-by-cyclic is virtually nilpotent, as witnessed by the fundamental group of a closed $3$-manifold admitting a $\textup{Sol}$ geometry~\cite[Corollary~12.7.7]{martelli}. However, one can still obtain the estimate $\cat_{\Poly}(M_f) \leq 4$. Indeed, if $M$ is flat (i.e., it admits the geometry~$\mathbb{E}^3$), then the fundamental group of $M$ is finitely generated virtually Abelian by a celebrated result of Bieberbach~\cite[Corollary~4.4.11]{martelli}. Similarly, if $M$ admits the $\textup{Nil}$ geometry, then $\pi_1(M)$ is finitely generated virtually nilpotent~\cite[Corollary~12.5.8]{martelli}. In both cases $\cat_{\Poly}(M) = 1$. Hence, applying Theorem~\ref{thm:mapping:tori:classical:estimate:LM} we get
    \[
    \cat_{\Poly}(M_f) \leq 2 \cdot \cat_{\Poly}(M) = 2 \cdot 1 = 2.
    \]
    We are left to consider the case in which $M$ admits the $\textup{Sol}$ geometry. In this case $\pi_1(M)$ is virtually solvable but \emph{not} virtually nilpotent~\cite[Corollary~12.7.7]{martelli}. Hence, we have to first estimate $\cat_{\Poly}(M)$. Recall that $M$ is either a torus bundle over the circle or it is obtained by glueing via an Anosov map two copies of the twisted $I$-bundle over the Klein bottle~\cite[Theorem~1.18]{aschenbrenner20123}. In the first case we have $\cat_{\Poly}(M) = 2$ because of Theorem~\ref{thm:mapping:tori:classical:estimate:LM}. In the second case, since the Klein bottle has finitely generated virtually Abelian fundamental group and we are glueing the two copies along a torus, by Proposition~\ref{prop:glueings:amcat:standard} we also have $\cat_{\Poly}(M) = 2$. So applying once more Theorem~\ref{thm:mapping:tori:classical:estimate:LM} to $M_f$ we obtain
    \[
    1 < \cat_{\Poly}(M_f) \leq 4. \qedhere
    \]
\end{proof}

\subsection{Mapping tori over hyperbolic prime smooth $3$-manifolds} Let $M$ be an oriented closed connected smooth $3$-manifold admitting a $\mathbb{H}^3$ geometry (i.e., a hyperbolic metric) and let $f\colon M\to M$ be an orientation-preserving self-diffeomorphism. To estimate the $\calG$-category of $M_f$  we will prove that it always admits a smooth regular circle foliation with finite holonomy groups. This is also true in the compact case and it is a consequence of the following result.
\begin{prop}[periodic self-maps and circle foliations]\label{prop:periodic:foliation}
    Let $M$ be an orientable connected compact smooth manifold with (possibly empty) boundary and let $f\colon M\to M$ be a \emph{periodic} self-diffeomorphism of pairs. Then $M_f$ admits a smooth regular circle foliation $\calF$ with finite holonomy groups such that $\partial M_f$ is vertical.
\end{prop}
\begin{proof}
Let $p \geq 1$ be the integer corresponding to the period of $f$ and let $\calF=\calF_f$ be the foliation by suspension defined in Example \ref{ex:foliation:suspension}. By construction, if $\partial M \neq \emptyset$, then the boundary $\partial M_f$ is vertical. 
Since $f$ is periodic of period~$p$, for every $(x,t) \in M\times\R$ we have that  $(x,t)\sim (x,t+p)$. Hence, the image of the leaf $\{x\}\times \R\subset M\times \R$ in the quotient~$M_f$ is a circle. Let $p_x$ be the period of $x$, i.e., the minimum integer~$k \geq 1$ such that $f^k(x)=x$, and set $q_x=p/p_x \in \Z$ (note that $p_x$ always divides $p$).  We claim that $\Hol(L_{[(x,0)]})$ is isomorphic to a quotient of $ \Z/q_x\Z$. Since all leaves are circles, we only need to compute the holonomy for a generator of the fundamental group, that is the image of the path $\{x\}\times [0, p_x]$. Let $\alpha$ be such a curve and let $T_0$ be the image of $M\times \{0\}$, which is transversal to the foliation~$\calF$.
Up to identifying the image of $M\times \{0\}$ with $M$, it follows from the definition of holonomy that  
\[
\hol^{T_0}(\alpha)= \textup{germ}_x(f^{p_x})\colon \textup{Diff}_x(M)\to\textup{Diff}_x(M). 
\]
Since $f^p=\id_M$, we have that $(\hol^{T_0}(\alpha))^{q_x}$ is the (germ of the) identity, concluding the proof. 
\end{proof}
By a classical application of Mostow's rigidity theorem, every self-\hspace{0pt}diffeomorphism of a complete finite-volume hyperbolic $n$-manifold with $n \geq 3$ is isotopic to a periodic one~\cite[Corollary~13.3.7]{martelli}.
Hence, we have the following corollary:
\begin{cor}[circle foliations of mapping tori over hyperbolic $3$-manifolds]\label{cor:hyperbolic:foliation}
    Let $M$ be an oriented compact connected smooth $3$-\hspace{0pt}manifold with (possibly empty) boundary, whose interior admits a complete finite-volume hyperbolic metric, and let $f\colon (M, \partial M) \to (M, \partial M)$ be an orientation-preserving self-diffeomorphism of pairs. Then, $M_f$ admits a smooth regular circle foliation $\calF$ with finite holonomy such that $\partial M_f$ is vertical.
\end{cor}
\begin{proof}
    The map $f$ is isotopic to a periodic self-diffeomorphism. Since isotopic maps give diffeomorphic mapping tori (Remark~\ref{rem:isotopic:map:give:isotopic:mapping:tori}), it is sufficient to apply Proposition~\ref{prop:periodic:foliation} .
\end{proof}

Finally, we deduce an estimate for the $\calG$-category:

\begin{cor}[$\cat_{\calG}$-estimate for mapping tori over hyperbolic $3$-manifolds]\label{cor:hyperbolic:cat:minent}
    Let $M$ be an oriented compact connected hyperbolic $3$-manifold with (possibly empty) boundary and let $f\colon (M, \partial M) \to (M, \partial M)$ be an orientation-preserving self-diffeomorphism of pairs. Then, we have
    \[ 1 < \cat_{\Am} (M_f) \leq \cat_{\Poly} (M_f)\leq 4.\]
    \
\end{cor}
\begin{proof}
    The fundamental group $\pi_1(M_f)$ is not amenable since it contains a subgroup isomorphic to $\pi_1(M)$, that is not amenable. Hence, we have \[\amcat(M_f) >1.\]
    On the other hand, applying Corollary~\ref{cor:hyperbolic:foliation} and Proposition~\ref{prop:amcat:circle:foliation:with:boundary} we immediately have \[\cat_{\Poly}(M_f)\leq 4. \qedhere\] 
\end{proof}
\begin{rem}[higher dimensions]\label{rem:hyp:higher:dim}
    Since the Mostow's rigidity argument holds in every dimension $d\geq3$, Corollary~\ref{cor:hyperbolic:foliation} and Corollary~\ref{cor:hyperbolic:cat:minent} (replacing $4$ with $d+1$) hold for all mapping tori over closed hyperbolic $d$-dimensional manifolds.
\end{rem}

\subsection{Mapping tori over aspherical Seifert fibered prime smooth $3$-mani\-folds}
Let $M$ be an oriented closed connected aspherical smooth $3$-manifold admitting a $\widetilde{\textup{SL}_2}$ or a $\mathbb{H}^2 \times \mathbb{R}$ geometry. Then $M$ admits a Seifert fibration~\cite[Theorem~12.8.1]{martelli}, i.e., a smooth regular circle foliation~$\calF$ with finite holonomy groups having a finite number of exceptional leaves. We will show that for every orientation-preserving self-diffeomorphism $f\colon M\to M$, the mapping torus  $M_f$ also admits a circle foliation with finite holonomy groups. This is also true in the compact case and is a corollary of the following proposition.
\begin{prop}[self-maps preserving the foliation induces a foliation in the mapping torus]\label{prop:fiber:preserving:MT:foliation}
    Let $M$ be an oriented compact connected smooth manifold with (possibly empty) boundary and let $\calF$ be a smooth regular foliation of $M$ such that $\partial M$ is vertical. Suppose that $f\colon (M, \partial M) \to (M, \partial M)$ is an orientation-preserving self-diffeomorphism of pairs such that $f(L)\in \calF$ for every $L\in \calF$. Then the mapping torus $M_f$ admits a foliation $\calF'$ such that $\partial M_f$ is vertical, every leaf $L'\in \calF'$ is diffeomorphic to some $L\in \calF$ and $\Hol(L')\cong \Hol(L)$.
\end{prop}
\begin{proof}
    As in Example~\ref{ex:foliation:suspension}, we identify $M_f$ with a quotient of $M\times \R$ under the $\Z$-action $k \cdot (x, t) = (f^k(x), t-k)$ for all $k \in \Z$ and $(x, t) \in M \times \R$. Let $\calP=\{t\}_{t\in \R}$ be the point foliation of $\R$. Then we can consider the product foliation on $M\times \R$:
    \[
    \calF\times\calP=\{L\times \{t\}\}_{L\in \calF, \ t \in \R}.
    \]
Since $f$ sends leaves of $\mathcal{F}$ to leaves of $\mathcal{F}$, this foliation is invariant under the $\Z$-action, therefore it induces a foliation $\calF'$ to the quotient space $M_f$. For all $t\in \R$ and $\varepsilon<1/2$ the projection map is an embedding when restricted to $M\times (t-\varepsilon, t+\varepsilon)$, hence each leaf $L\times\{t\}$ of $\calF\times\calP$ is diffeomorphic to its image $L'$ in the mapping torus. To compute the holonomy of $L'$, we can take a transverse submanifold contained in the image of $M\times (t-\varepsilon, t+\varepsilon)$. It is then straightforward to check that
\[
\Hol(L')\cong \Hol(L\times\{t\}) \cong \Hol(L).
\]   
By construction, if $\partial M$ is vertical with respect to $\calF$, also  $\partial M\times \R$ is vertical with respect to $\calF\times\calP$, hence $\partial M_f$ is vertical with respect to $\calF'$.
\end{proof}

The previous proposition allows us to prove that mapping tori over specific Seifert fibered smooth $3$-manifolds admit a smooth regular circle foliation:

\begin{cor}[mapping tori over Seifert fibered $3$-manifolds admit circle foliations]\label{cor:SF:MT:foliation}
   Let $M$ be an oriented compact Seifert fibered manifold with (possibly empty) boundary which falls in one of the following cases:
   \begin{itemize}
       \item $M$ is closed and it admits either a $\widetilde{\textup{SL}_2}$ or a $\mathbb{H}^2 \times \mathbb{R}$ geometry;
       \item $M$ has non-empty boundary and is not diffeomorphic to $S^1 \times D^2$ or $T^2\times I$.
       \end{itemize}
   Then, for every orientation-preserving self-diffeomorphism $f \colon (M, \partial M) \to (M, \partial M)$ of pairs, the mapping torus~$M_f$ admits a smooth regular circle foliation with finite holonomy groups such that the (possibly empty) boundary~$\partial M_f$ is vertical.
\end{cor}
Note that this statement also includes the ``twisted'' interval-bundle over the Klein bottle $M =K\tilde{\times}I$, which is a piece that can occur in JSJ decompositions which will be discussed in the next subsection (Subsection~\ref{subsec:JSJ}).
\begin{proof}
    Under the hypothesis of the corollary, the Seifert fibered manifold $M$ has the following property: Every orientation-preserving self-diffeomorphism $f \colon M \to M$ is isotopic to a fiber-preserving orientation-preserving self-\hspace{0pt}diffeomorphism
    \cite[Theorem 3.11]{jiang1996homeomorphisms}.
    Hence, since isotopic maps give diffeomorphic mapping tori (Remark~\ref{rem:isotopic:map:give:isotopic:mapping:tori}), Proposition~\ref{prop:fiber:preserving:MT:foliation} concludes the proof.
\end{proof}


Finally, we can estimate the $\calG$-category for mapping tori over oriented closed connected aspherical smooth $3$-manifolds admitting a $\widetilde{\textup{SL}_2}$ or a $\mathbb{H}^2 \times \mathbb{R}$ geometry:
\begin{cor}[$\cat_{\calG}$-estimate for mapping tori over aspherical Seifert fibred closed $3$-manifolds]\label{cor:MT:aspherical:Seifert}
    Let $M$ be an oriented closed connected $3$-manifold admitting a $\mathbb{H}^2\times\R$ or a $\widetilde{\textup{SL}_2}$ geometry and let $f\colon M\to M$ be an orientation-preserving self-diffeomorphism. Then, we have
    \[ 1< \amcat (M_f) \leq \cat_{\Poly} (M_f)\leq 4.\]
\end{cor}
\begin{proof}
    Since the fundamental group of $M$ is either virtually isomorphic to $\Z\times F$ or it fits in the short exact sequence
    \[
    1 \to \Z \to \pi_1(M)\to F\to 1
    \] where $F$ is some non-Abelian free group~\cite[Table 1]{aschenbrenner20123}, it is not amenable. Then also $\pi_1(M_f)$ is not amenable, hence
    \[\amcat(M_f)>1.\]
    On the other hand, applying Corollary~\ref{cor:SF:MT:foliation} and Proposition~\ref{prop:amcat:circle:foliation:with:boundary} we immediately have
    \[\cat_{\Poly}(M_f)\leq 4. \qedhere\]
\end{proof}

\subsection{Mapping tori over irreducible smooth $3$-manifolds with non-trivial JSJ-decomposition}\label{subsec:JSJ}
We are left to consider the case where $M$
is an oriented closed connected irreducible smooth $3$-manifold that does not admit a Thurston geometry. Let $f\colon M\to M$ be an orientation-preserving self-diffeomorphism of $M$. Recall that an embedded torus $T^2 \subset M$ is said to be \emph{incompressible} if the morphism induced by the inclusion $\pi_1(T^2) \to \pi_1(M)$ is injective. Moreover, a compact manifold~$N$ with toroidal boundary is said to be \emph{atoroidal} if every map from a torus $T^2$ to $N$ can be homotoped to the boundary~$\partial N$. By a classical result by Jaco, Shalen, and Johansson~\cite[Theorem~1.6.1]{aschenbrenner20123}, $M$ admits a \emph{JSJ-decomposition}:
There exists a minimal collection of disjoint embedded incompressible tori $\mathcal{T}=\{T_1,\cdots,T_k\}$
such that each piece obtained by cutting along the tori (i.e., each connected component of $M \cut \mathcal{T}$) is either atoroidal  or Seifert fibered. This collection of tori is unique up to isotopy.
Furthermore, we can assume that each atoroidal piece is hyperbolic (i.e., its interior admits a hyperbolic metric of finite volume), since the Hyperbolization Theorem~\cite[Theorem~1.7.5]{aschenbrenner20123} ensures that such pieces are either hyperbolic or with finite fundamental group (and we can exclude the latter case, since this would imply that the piece is itself already closed, whence the manifold is elliptic~\cite[Theorem~1.7.3]{aschenbrenner20123}). Finally, since we are assuming that $M$ does not admit a Thurston geometry, $M$ cannot be a torus bundle~\cite[Theorem~1.10.1]{aschenbrenner20123} and so all the Seifert fibered pieces appearing in the JSJ-decomposition of~$M$ satisfy the hypothesis of Corollary~\ref{cor:SF:MT:foliation}~\cite[Proposition~1.6.2]{aschenbrenner20123}.

We are now ready to study in detail the topology of mapping tori over non-geometric irreducible manifolds. The uniqueness up to isotopy of $\mathcal{T}$ implies that, up to replacing $f$ with an isotopic orientation-preserving self-diffeomorphism, we can assume that $f(\mathcal{T})= \mathcal{T}$. Since mapping tori of isotopic maps are diffeomorphic (Remark~\ref{rem:isotopic:map:give:isotopic:mapping:tori}), we assume we are in the situation described in the following setup:
\begin{setup}\label{setup:JSJ}
We consider the following situation:
\begin{itemize}
    \item $M$ is an oriented closed connected irreducible smooth $3$-manifold that does not admit a Thurston geometry;
    \item $\mathcal{T} = \{T_1, \cdots, T_k\}$ is a minimal collection of disjoint embedded incompressible tori in the (unique up to isotopy) JSJ decomposition of $M$;
    \item $f \colon M \to M$ is an orientation-preserving self-diffeomorphism such that $f(\mathcal{T}) = \mathcal{T}$;
    \item Denote by $T=\cup_{i=1}^kT_i$ the union of all the tori appearing in the collection~$\mathcal{T}$ and by $f_T$ the restriction of $f$ to $T$;
    \item Then the mapping torus $S:=T_{f_T}$ is an embedded closed (possibly disconnected) codimension-one submanifold of $M_f$. Now we have that $S\cong \sqcup_{i=1}^s S_i$ for some $s\leq k$, where each $S_i$ is diffeomorphic to the mapping torus of some self-diffeomorphism of the $2$-torus (Corollary~\ref{cor:MT:disconnected:general});
    \item Let $\mathcal{S} =\{S_1, \cdots, S_s\}$ be the collection of the connected components of $S$.
    \end{itemize}
\end{setup}

\begin{lem}[decomposition of $M_f$ in connected components]\label{lem:mapping:tori:decomposition:JSJ}
    In the situation of Setup~\ref{setup:JSJ} we have that every connected component of $M_f \cut \mathcal{S}$ is diffeomorphic to the mapping torus of an orientation-preserving self-diffeomorphism of a connected component of $M\cut\mathcal{T}$.
\end{lem}
\begin{proof}
First, we can extend by continuity the restriction of $f$ to $M \setminus T$ (that is diffeomorphic to the interior of $M\cut \mathcal{T})$  to an orientation-preserving self-diffeomorphism $\hat{f}\colon M \cut \mathcal{T} \to M \cut \mathcal{T}$.
Secondly, since  $M_f \backslash S$ is exactly the mapping torus of $f |_{M\setminus T}$, we have that $M_f {\cut}  {\mathcal{S}}$ 
is diffeomorphic to $(M \cut \mathcal{T})_{\hat{f}}$. 
Thus, we have that each connected component of $M_f\cut\mathcal{S}$ is diffeomorphic to the mapping torus of an orientation-preserving self-diffeomorphism $g_N\colon N\to N$, where $N$ is a connected component of $M\cut \mathcal{T}$ (Corollary~\ref{cor:MT:disconnected:general}). More precisely $g_N$ is the composition of the restrictions of $\hat{f}$ to the connected components in the orbit of $N$. 
\end{proof}

We can use the decomposition in Lemma~\ref{lem:mapping:tori:decomposition:JSJ} to estimate the $\calG$-category of the mapping torus $M_f$:
\begin{prop}[$\cat_{\calG}$-estimate for mapping tori over irreducible $3$-mani\-folds with non-trivial JSJ-decomposition]\label{prop:MT:JSJ:catam:minent}
    Let $M$ be an oriented closed connected irreducible smooth $3$-manifold that does not admit a Thurston geometry, and let $f\colon M\to M$ be an orientation-preserving self-diffeomorphism. Then, we have that
    \[ \amcat(M_f)\leq \cat_{\Poly}(M_f)\leq 4.\]
\end{prop}
\begin{proof}
    In the situation of Setup~\ref{setup:JSJ}, Lemma~\ref{lem:mapping:tori:decomposition:JSJ} shows that
    \[
    M_f\cut \mathcal{S}\cong  (N_1)_{g_1}\sqcup \dots \sqcup (N_s)_{g_s},
    \]
    where for every $i \in \{1,\cdots,s\}$ the manifold $N_i$ is a connected component of $M\cut \mathcal{T}$ and $g_i$ is an orientation-preserving self-diffeomorphism. Moreover, as discussed above, each $N_i$ is either a hyperbolic manifold or a Seifert fibered manifold satisfying the hypothesis of Corollary~\ref{cor:SF:MT:foliation}. Then, by Corollary~\ref{cor:hyperbolic:foliation} and Corollary~\ref{cor:SF:MT:foliation}, respectively, we have that for every $i\in\{1,\cdots,s\}$ the mapping torus $(N_i)_{g_i}$ admits a smooth regular circle foliation such that the boundary is vertical.
    Since each component of the boundary $\partial (N_i)_{g_i}$ is diffeomorphic to a mapping torus of some orientation-preserving diffeomorphism of the $2$-torus, every component of $\partial (N_i)_{g_i}$ admits an $\mathbb{E}^3$, a $\textup{Nil}$ or a $\textup{Sol}$ geometry~\cite[Theorem~1.10.1]{aschenbrenner20123}.
    In the case where $N_i$ is hyperbolic, the map $g_i$ is isotopic to a periodic map. This implies that $(N_i)_{g_i}$ is finitely covered by the product~$N_i \times S^1$ and thus the boundary components~$\partial(N_i)_{g_i}$ are finitely covered by the disjoint union of some copies of a $3$-torus, i.e., each of them admits a $\mathbb{E}^3$ geometry. On the other hand, in the case where $N_i$ is Seifert fibered, we have that the smooth regular circle foliation of $\partial(N_i)_{g_i}$ has finite holonomy groups, since they are isomorphic to the holonomy groups of $\partial N_i$ by  Corollary~\ref{cor:SF:MT:foliation}.
    In other words $\partial(N_i)_{g_i}$ is Seifert fibered, hence it does not admit a $\textup{Sol}$ geometry~\cite[Theorem~12.1.1]{martelli}.
    
    Since closed $3$-manifolds admitting either a $\mathbb{E}^3$ or a $\textup{Nil}$ geometry have finitely generated virtually nilpotent fundamental group, we have that $\partial (N_i)_{g_i}$ is a $\Poly$-set in $(N_i)_{g_i}$. Thus, as discussed in Section \ref{Sec:amcat:glueings}, the submanifolds $S_i$ are either diffeomorphic to or doubly covered by  some  $\partial(N_i)_{g_i}$, hence their fundamental group still lies in $\Poly$.
    By Corollary~\ref{cor:final:glueing:foliations} we conclude that \[\cat_{\Am}(M_f) \leq \cat_{\Poly}(M_f)\leq 4. \qedhere \]
    %
\end{proof}

\subsection{Proof of Theorem~\ref{thm:amcat:main}}
It is sufficient to combine Proposition~\ref{prop:MT:non:aspherical}, Proposition~\ref{prop:MT:infinite:amenable}, Corollary~\ref{cor:hyperbolic:cat:minent}, Corollary~\ref{cor:MT:aspherical:Seifert} and Proposition~\ref{prop:MT:JSJ:catam:minent}.

\section{Minimal volume entropy of mapping tori over $3$-manifolds}\label{Sec:minent:MT}

Building on the computations of the minimal volume entropy of mapping tori over prime $3$-manifolds (Corollary~\ref{cor:minent:MT:prime}), in this section we show how to extend the vanishing result to all mapping tori over $3$-manifolds.

\begin{thm}[$\minent$ of mapping tori over $3$-manifolds]\label{thm:minent_reducible}
     Let $M$ be an oriented closed smooth $3$-manifold.
     Then, for every orientation-preserving self-diffeomorphism~$f \colon M \to M$, the mapping torus $M_f$ has zero minimal volume entropy: 
    \[
    \minent(M_f) = 0.
    \]
\end{thm}

In the proof of \Cref{thm:minent_reducible} we work with (a power of) the volume entropy of homology classes~\cite[p.~4395]{BS2023seminorm}. Following Babenko and Sabourau, given a path-connected topological space~$X$, we define the \emph{volume entropy of a homology class} $\alpha \in H_n(X; \Z)$ as
\[
\omega(\alpha) := \inf_{(P, f)} \omega_{\ker f_*}(P) \in [0, \infty),
\]
where the infimum is taken over all orientable closed $n$-pseudomanifolds~$P$, and all continuous maps $f \colon P \to X$ such that $H_n(f)([P]) = \alpha$ (essentially, the infimum is taken over the \emph{geometric cycles} representing $\alpha$~\cite[p.~108--109]{hatcher}).
Here $\omega_{\ker f_*}(P)$ denotes the minimal volume entropy of $P$ associated to the normal subgroup $\ker(f_*) \trianglelefteq \pi_1(P)$, that is, if $P_H$ is the normal covering of $P$ associated to $H \trianglelefteq \pi_1(P)$ we set 
\[
\omega_{H}(P) = \inf_{g} \ent_H(P, g) \cdot \vol(P, g)^{\frac{1}{n}},
\]
where $\ent_H(P, g) = \lim_{R \to +\infty} (1 \slash R) \log (\vol_g(B(R)))$ is the entropy of the normal covering $P_H$ (namely, the balls lie in $P_H$) and $g$ runs over all the piecewise-Riemannian metrics on $P$, i.e., a family of Riemannian
metrics on all simplices of~$P$ that is compatible along common subsimplices~\cite[Definition~6.1]{lmfibration}.

Clearly, when $P$ is a manifold and $H$ is the trivial group, we have $\omega_{H}(P) = \minent(P)$. This readily shows that 
\[
\omega([M]) \leq \minent(M)
\]
for all oriented closed connected smooth $n$-manifolds~$M$. It is often convenient to work with the $n$-th power of $\omega$ instead of $\omega$ itself. Indeed, as proved by Babenko and Sabourau, the real function
\begin{align*}
 \Omega \colon H_n(X;\mathbb{Z})&\to [0,+\infty)   \\
 \alpha &\mapsto \omega(\alpha)^n
\end{align*}
gives a pseudo-distance on $H_n(X; \mathbb{Z})$ and satisfies the following properties:

\begin{prop}[properties of $\Omega$]\label{prop:properties:omega}
The function $\Omega$ satisfies the following properties:
\begin{enumerate}
    \item Given a path-connected space~$X$ and two homology classes $\alpha, \beta \in H_n(X; \mathbb{Z})$, we have \[\Omega(\alpha + \beta) \leq \Omega(\alpha) + \Omega(\beta);\]
    \item Given a continuous map~$f \colon X \to Y$ between path-connected topological spaces and a homology class $\alpha \in H_n(X; \Z)$, we have \[\Omega\big(H_n(f)(\alpha)\big) \leq \Omega(\alpha);\]
    \item For every oriented closed connected smooth $n$-manifold 
    we have \[\minent(M)^n = \Omega\big([M]\big),\] where $[M]\in H_n(M;\mathbb{Z})$ is the fundamental class of $M$;
    \item\label{it:omega_classifying} For every oriented closed connected smooth $n$-manifold  
    we have
    \[\minent(M)^n = \Omega\big(H_n(c_M)([M])\big),\]
    where $c_M \colon M \to K(\pi_1(M), 1)$ denotes the classifying map of $M$;
    \item Let $M$ be an oriented closed connected $n$-manifold and let $[M] \in H_n(M;\mathbb{R})$ be its fundamental class.
    Let $X$ be a connected CW-complex and suppose that there exists a continuous map $f \colon M \to X$ inducing an isomorphism between the fundamental groups.
    Then, we have
    \[\Omega\big(H_n(f)([M])\big) = \Omega([M]) = \minent(M)^n.\]
\end{enumerate}
\end{prop}
\begin{proof}
    \emph{Ad~1 \& 2)} These facts were proved by Babenko and Saborau~\cite[Equation~(1.3) and Theorem~1.1]{BS2023seminorm};

    \emph{Ad~3)} This follows from the fact that if a geometric cycle is represented by a manifold and the associated map is $\pi_1$-surjective, then that geometric cycle realizes the infimum~\cite[p.~1114]{brunnbauer2008homological};

    \emph{Ad~4)} The two inequalities can be deduced from item~(\emph{2}) and a result by Brunnbauer~\cite[Theorem~10.2]{brunnbauer2008homological};

    \emph{Ad~5)} By items~(\emph{3}) and~(\emph{4}) we have \[\minent(M)^n = \Omega([M]) = \Omega(H_n(c_M)([M])).\] 
    Hence, it is sufficient to consider the induced map $\varphi \colon X \to K(\pi_1(M), 1)$ corresponding to the group homomorphism~$\varphi_* = c_M^* \circ f_*^{-1}$. The comparison map~$c_M$ of $M$ is homotopic to the composition $\varphi \circ f$. Item~(\emph{2}) concludes the proof.
\end{proof}

Item~(\emph{3}) of the previous proposition shows that the vanishing of the minimal volume entropy of an oriented closed connected smooth manifold~$M$ 
is equivalent to $\Omega([M]) = 0$.
Properties~(\emph{3}) and~(\emph{4}) have been often stated in the literature with the assumption that the dimension is $n \ge 3$. However, they hold also for $n < 3$, as it follows readily from a result of Babenko and Sabourau~\cite[Theorem 3.1 (6)]{BS2023seminorm}; in any case, we will need these properties for $n = 4$ only.

\subsection{Diffeomorphisms of 3-manifolds}\label{sec:homeo3}
Before going through the proof of \Cref{thm:minent_reducible}, we need to recall a classical result by McCullough~\cite[Section~3]{McCull} that shows that the mapping class group of an oriented closed connected smooth $3$-manifold is generated by four types of mapping classes. We work with the following setup:

\begin{setup}\label{setup:reducible:manifolds}
    Let $M$ be an oriented closed connected smooth $3$-manifold. By the Kneser-Milnor theorem~\cite{Kneser, Milnor}, there exist two natural numbers $k$ and $h$ such that $M$ can be written as a connected sum
\[M = M_1 \connsum \cdots \connsum M_k \connsum \big(\connsum_h S^2 \times S^1\big),\]
where, for each $i\in\{1, \ldots, k\}$,  $M_i$ is an oriented closed connected irreducible smooth $3$-manifold. Following the survey paper of McCullough~\cite[Section~3]{McCull}, we view $M$ as a punctured 3-sphere $W$ with $k+2h$ spherical boundary components 
\[S_1, \ldots, S_k, S_{k+1,0}, S_{k+1,1}, \ldots, S_{k+h,0}, S_{k+h,1}\]
on which the pieces of the prime decomposition of $M$ are attached as follows:
\begin{itemize}
    \item For each $i\in\{1, \ldots, k\}$ we pick a $3$-ball $D_i \subset M_i$, and we attach $M_i':=M_i \setminus\interior(D_i)$ along $\partial D_i$ to the sphere $S_i \subset\partial W$;
    \item For each $j\in\{1, \ldots, h\}$ we attach a copy of $S^2 \times [0,1]$ to $W$ by identifying $S^2 \times \{0\}$ with $S_{k+j,0}$ and $S^2\times\{1\}$ with $S_{k+j,1}$.
\end{itemize}
Henceforth, we interpret $W$, the $M_i'$'s and the copies of $S^2\times [0,1]$ as subsets of $M$.
In order to distinguish the different copies of $S^2\times [0,1]$ we denote them by $S_j^2 \times [0,1]$ with $j \in \{1,\cdots,h\}$.
    Moreover, we fix the following objects and notation:
    \begin{itemize}
        \item For every $i \in \{1,\cdots, k\}$ we denote by $C_i$ a fixed collar of $S_i = \partial M_i'$ inside $M_i'$.
        By construction $C_i \cong S^2 \times [0, 1]$ and we fix an appropriate diffeomorphism that identifies $S^2 \times \{1\} = S_i \subset \partial W$;

        \item For every $j \in \{1, \cdots, h\}$ we fix two disjoint collars~$C_{j,0}$ and $C_{j,1}$ of $S_{j, 0}$ and $S_{j, 1}$ inside $S_j^2 \times [0,1]$. Similarly as before, we identify such collars with two copies of $S^2 \times [0, 1]$ such that the first diffeomorphism gives $S^2 \times \{1\} = S_{k+j, 0}$ and the second one gives $S^2 \times \{1\} = S_{k+j, 1}$. 
    \end{itemize}
        Since all these collars lie in $M \setminus \textup{int}(W)$, it is readily seen that they are all pairwise disjoint. Moreover, they all intersect $W$ at the boundary $\partial W$ only.
\end{setup}

In the situation of Setup~\ref{setup:reducible:manifolds}, we can describe four types of \emph{elementary} diffeomorphisms:
\begin{itemize}
    \item \textit{Type 1. Diffeomorphisms preserving summands.} The restriction of such a diffeomorphism to $W$ is the identity (thus, such a diffeomorphism is induced by self-diffeomorphisms of the pieces $M_i'$ and $S_j^2\times [0,1]$ fixing their boundaries);
    \item \textit{Type 2. Interchange of diffeomorphic summands.} These diffeomorphisms interchange two pieces $M'_i$ and $M'_j$ which are diffeomorphic (or two pieces $S_{j_1}^2 \times [0,1]$ and $S_{j_2}^2\times [0,1]$) via an orientation-preserving 
    diffeomorphism and keep all the other summands fixed, leaving $W$ invariant.
    In the case where two pieces $S_{j_1}^2 \times [0,1]$ and $S_{j_2}^2\times [0,1]$ are exchanged, this is done via the natural map that reads as the identity map on $S^2\times [0,1]$;
    \item \textit{Type 3. Spins of $S^1 \times S^2$ summands.}
    Let $j \in \{1, \ldots, h\}$.
    These diffeomorphisms interchange $S_{k+j,0} = S^2_j\times \{0\}$ with $S_{k+j,1} = S^2_j\times\{1\}$, i.e., they restrict to an orientation-preserving diffeomorphism that interchanges the boundary components of $S_j^2 \times [0,1]$, fixing all the other summands and leaving $W$ invariant;
    \item \textit{Type 4. Slide diffeomorphisms.}
    Let $S$ be one of the boundary spheres of $W$, which either bounds $M_i'$ for some $i \in \{1,\cdots,k\}$ or is one of the boundary components of $S_j^2\times [0,1]$ for some $j \in \{1,\cdots,h\}$.
    Let $\alpha\colon [0,1]\to M$ be an \emph{embedded smooth} arc with endpoints lying in $S$ such that \[\alpha((0, 1)) \subset M \setminus M_i' \quad \big(\text{resp.} \, \, \alpha((0, 1)) \subset M \setminus (S_j^2 \times [0,1])\big)\] if $S = \partial M_i'$ (resp.\ if $S \subset \partial (S_j^2 \times [0,1])$).
    Let $N \subset N'$ be two nested closed tubular neighborhoods of the union $\alpha([0,1]) \cup S$ such that $N \subset \interior (N')$.
    By construction the closed set~$N' \setminus \interior(N)$ consists of two connected components: A component is diffeomorphic to $S^2 \times [0,1]$ and the other is diffeomorphic to $T^2 \times [0,1]$. Under these identifications we have that $S^2\times\{0\} \cup T^2\times\{0\} = \partial N$ and $S^2\times\{1\} \cup T^2\times\{1\} = \partial N'$.
    We set $T := T^2 \times [0,1]$.
    Let $\hat\alpha$ be a closed curve obtained by closing up $\alpha$ with an arc that is entirely contained in $S$. Let $S^1_{\hat\alpha}$ be a curve in $T$ that is freely homotopic to $\hat{\alpha}$ in $N'$. Using this new curve~$S^1_{\hat\alpha}$ we can make the following identification: $T \cong (S^1_{\hat\alpha} \times S^1) \times [0,1]$. 

    A \emph{slide diffeomorphism} associated to $\alpha$ is the identity on $M \setminus T$ and on $T$ it is defined by $(x,y,z)\mapsto (x+z,y,z)$, where we are identifying $S^1_{\hat\alpha}$ with $\mathbb{R}/\mathbb{Z}$. 
\end{itemize}

\begin{rem}\label{rem:regular:neigh:not:important}
    The (isotopy class of the) resulting slide diffeomorphism does not depend on the choice of the nested tubular neighbourhoods~$N$ and~$N'$.
    In the same way, the arc $\alpha$ only matters up to isotopy. 
\end{rem}

\begin{rem}\label{rem:elementary_collar}
    For diffeomorphisms of types 1, 2 and 3, which preserve $W$, we further assume that they preserve the collars $C_i$, $C_{j,0}$, $C_{j,1}$ possibly exchanging them, but preserving the fixed identifications with $S^2\times [0,1]$ (i.e., up to these identifications, the restrictions to the fixed collars all read as identity maps of $S^2\times [0,1]$).
    We can safely do so, up to changing the elementary diffeomorphism with an isotopic one: this follows from the well-known facts that orientation-preserving diffeomorphisms of $S^2$ are isotopic to the identity, and uniqueness of collars up to isotopy.
\end{rem}

The following result is due to McCullough.

\begin{thm}[{\cite[p.~69]{McCull},~\cite[Theorem 4.1]{bucherneofytidis},~\cite[Section~3.1.1]{mann2020dynamical}}]\label{thm:mccullough}
    Let $M$ be an oriented closed connected smooth $3$-manifold, and let $f\colon M \to M$ be an orientation-preserving diffeomorphism.
    Then, $f$ is isotopic to a composition
    \[g_4 \circ g_3 \circ g_2 \circ g_1,\]
    where, for each $\ell\in \{1,2,3,4\}$, the map $g_\ell$ is an orientation-preserving self-diffeomorphism, equal to the composition of finitely many diffeomorphisms of type $\ell$.
\end{thm}

\subsection{Proof of Theorem \ref{thm:minent_reducible}}
Throughout the proof, we consider the situation of Setup~\ref{setup:reducible:manifolds}. Given the prime decomposition
\[
    M = M_1 \connsum \cdots \connsum M_k \connsum \big(\connsum_h S^2 \times S^1\big),
\]
it is convenient to identify $M_i$ with $M_i' \slash \partial M_i'$ (i.e., we are identifying the boundary of $M_i'$ to a point). 

In the proof, we use some auxiliary CW-complexes associated to $M$, which we now introduce.
    Starting from the prime decomposition of $M$, we define the following CW-complex:
    \[M^\vee:=M_1 \vee \ldots \vee M_k \vee \bigvee_h S^1.\]
    Recall, from \Cref{setup:reducible:manifolds}, that $M$ is decomposed into a punctured sphere $W$ and pieces $M_i'$, $S^2_j\times[0,1]$; 
    then, the CW-complex~$M^\vee$ is obtained from $M$ by first collapsing the punctured sphere~$W \subset M$  to a point and secondly collapsing in each $S_j^2 \times [0,1]$ piece the $S^2$ factor to a point. We denote by $\pi \colon M \to M^\vee$ the continuous map that performs such collapses.
    We also denote for convenience the wedge-sum basepoint of $M^\vee$ by $x_0$.
    By construction, if $w$ is a baseboint of $M$ lying in $W$, the quotient map $\pi$ induces an isomorphism between the fundamental groups $\pi_1(M, w)$ and $\pi_1(M^\vee, x_0)$.

    The other complex (in fact, a manifold, possibly not connected) we need to consider is the disjoint union of $M_1,\dots,M_k$, which we denote by $M^\sqcup$.
    We have the map $\iota \colon M^\sqcup \to M^\vee$ that sends each $M_i$ to its corresponding copy in $M^\vee$ (mapping all the basepoints to $x_0 \in M^\vee$).

    On $M^\vee$ and $M^\sqcup$ we fix some preferred homology classes, which we call \emph{fundamental classes}.
    For $M^\vee$, we consider 
    \[[M^\vee]:=\sum_{i = 1}^k H_3(\iota_i)([M_i]) \in H_3(M^\vee;\mathbb{Z}),\]
    where $\iota_i \colon M_i \to M^\vee$ denotes the inclusion of the $i$-th $M_i$ factor into $M^\vee$. A straightforward computation of local degrees shows that $H_3(\pi)([M]) = [M^\vee]$ which justifies the terminology of fundamental class of $M^\vee$.
    
    On $M^\sqcup$, we consider $[M^\sqcup] := \sum_{i = 1}^k [M_i]\in H_3(M^\sqcup;\mathbb{Z})$; by construction, we have $H_3(\iota)([M^\sqcup]) = [M^\vee]$.
    



We divide the proof of \Cref{thm:minent_reducible} in two steps, in which we describe how a self-diffeomorphism $f\colon M \to M$ induces certain self-maps $f^\vee$ and $f^\sqcup$ of $M^\vee$ and $M^\sqcup$, respectively, so that we have maps between the corresponding mapping tori, as in this diagram:
\begin{equation}\label{eq:schema}  M_{f} \rightarrow M^\vee_{f^\vee} \leftarrow M^\sqcup_{f^{\sqcup}}. \end{equation}
Using the general machinery presented in \Cref{sec:prelim:MT}, we shall see that the fundamental classes of $M$, $M^\vee$ and $M^\sqcup$ lift to ``fundamental classes'' on the respective mapping tori, and that they correspond to each other under the maps in \eqref{eq:schema}.
Then, we shall conclude by using the properties of $\Omega$ listed in \Cref{prop:properties:omega}.



\medskip

        \textbf{(Step 1)\quad}
%
%
    Let $f \colon M \to M$ be a self-diffeomorphism.
    In this part of the proof, we construct a continuous self-map~$f^\vee \colon M^\vee \to M^\vee$ such that the diagram
    \begin{equation}\label{eq:diagram_f}\begin{tikzcd}
        M \arrow["f", r] \arrow["\pi", d] & M \arrow["\pi", d]\\
        M^\vee \arrow["f^\vee", r] & M^\vee
    \end{tikzcd}\end{equation}
    commutes up to homotopy;
    as explained in Subsection~\ref{subsection:digram:induces:map}, this gives a map between the corresponding mapping tori on $M$ and $M^\vee$.
    %
    
    By Theorem~\ref{thm:mccullough}, we can assume that $f$ is a finite composition of diffeomorphisms of the four elementary types described in \Cref{sec:homeo3}. We fix such a sequence of elementary diffeomorphisms whose composition gives~$f$.


        The map $\pi$ maps a collar~$C_i \subset M_i'$ to a $3$-dimensional ball inside $M_i = \pi(M_i')$ by collapsing $\partial M_i'$ to $x_0$ and it sends a collar $C_{j,\epsilon} \subset S_j^2 \times [0, 1]$ to a segment in $S^1 = \pi(S_j^2 \times [0, 1])$ whose final endpoint agrees with $x_0 \in M^\vee$.
    
    We are now ready to describe the construction of the map~$f^\vee$. 
    First, assume that $f$ is just a diffeomorphism of any of the four elementary types.
    \begin{itemize}
        \item \textit{Types 1, 2, and 3.}
        In these cases, $f$ leaves $W$ invariant, and this implies that there is a unique map $f^\vee$ such that Diagram~\eqref{eq:diagram_f} strictly commutes. We define $f^\vee$ to be such a map. By construction, the restriction of $f^\vee$ to each $M_i \subset M^\vee$ is a (pointed at $x_0$) 
        orientation-preserving diffeomorphism (see \Cref{rem:elementary_collar}) from $M_i$ to some (not necessarily distinct) diffeomorphic copy~$M_{i'}\subset M^\vee$.
        \item \textit{Type 4.} 
        Suppose that $f$ is a slide diffeomorphism along an arc $\alpha \colon [0,1] \to M$ with endpoints in a sphere component $S \subset \partial W$.
        Recall from Setup~\ref{setup:reducible:manifolds} that we have fixed a collar $C = S^2\times[0,1]$ of $S$ in some~$M_i'$ or $S^2_j\times[0,1]$.
        The map $f^\vee$ is then defined as follows:
        \begin{itemize}
            \item $f^\vee(p) = p$ for $p \not\in \pi(C)$;
            \item $f^\vee(\pi(x,t)) = \pi(x,2t)$ for $(x,t) \in S^2 \times [0,1/2] \subset C$;
            \item $f^\vee(\pi(x,t)) = \pi(\alpha(2t-1))$ for $(x,t) \in S^2 \times [1/2,1] \subset C$.
        \end{itemize}
        This definition is well posed and gives a continuous self-map of $M^\vee$.
    \end{itemize}
    Finally, if $f$ is a finite composition of elementary diffeomorphisms, $f^\vee$ is defined as the composition of the corresponding self-maps of $M^\vee$ just described.

    \begin{lem}\label{lemma:first:homotopy:square:commute}
        Diagram~\eqref{eq:diagram_f}
        commutes up to homotopy.
    \end{lem}
    \begin{proof}
        For the proof, we can assume that $f$ is a diffeomorphism of any of the four elementary types; the conclusion will then hold also for compositions of such maps since $\pi \circ f_1 \circ f_2 \simeq f^\vee_1 \circ \pi \circ f_2 \simeq f^\vee_1 \circ f_2^\vee \circ \pi$.
        For diffeomorphisms of type \emph{1}, \emph{2} or \emph{3}, the diagram strictly commutes by construction. Therefore, we are left to consider the case in which $f$ is a slide diffeomorphism. Let $\alpha$ be the arc defining the slide diffeomorphism~$f$.
        
        In the situation of Setup~\ref{setup:reducible:manifolds}, by Remark~\ref{rem:regular:neigh:not:important} we can assume that the tubular neighbourhood~$N'$ of $\alpha([0,1]) \cup S$ contains the collar $C$ of $S$ that was used to define the map $f^\vee$.
        Since the slide diffeomorphism~$f$ is the identity outside $N'$, by construction the maps $\pi\circ f$ and $f^\vee\circ \pi$ both coincide outside of $N'$ with the quotient map $\pi$.
        So we are left to show that $\pi\circ f$ and $f^\vee\circ \pi$ are homotopic inside $N'$ relatively to the boundary of $N'$. Such restrictions can be reinterpreted as maps $N' \to \pi(N')$.
    
        Recall that $N'$ is diffeomorphic to a solid torus from which an open $3$-\hspace{0pt}dimensional ball has been removed.
        Since the arc $\alpha$ is embedded, up to changing $\alpha$ and $N'$ by isotopies (with a transversality argument) we can assume that the connected components of $N' \cap \partial W$ are $S$ and, possibly, only other meridian discs of $N'$, cutting $N'$ in ``slices'' that are alternately contained in $W$ and in pieces $M_j'$ or $S_j^2\times [0,1]$.
        It follows from this description and the definition of $\pi$ that the image $\pi(N')\subset M^\vee$ is homotopy equivalent to a wedge of circles based at $x_0$. Hence $\pi(N')$ is an aspherical space.

        As discussed above, the restriction $\pi\circ f$ to $N'$ is a map $u_0 \colon N'\to \pi(N')$, and, similarly, the restriction of $f^\vee\circ\pi$ to $N'$ is a map $u_1 \colon N'\to\pi(N')$. We already noticed that these maps agree on $\partial N'$, and so we are left to prove that they are homotopic relative to $\partial N'$ (so that the homotopy can be extended over the whole $M$ to a homotopy between $\pi \circ f$ and $f^\vee \circ \pi$).
        
        We fix a cell decomposition of $N'$ as in \Cref{fig:cell} consisting of two $0$-cells, four $1$-cells, four $2$-cells and one $3$-cell.
        All the $1$-cells of $N'$, except the red one in the picture, are contained in $\partial N'$.
        \begin{figure}[ht]
            \centering
            \includegraphics[scale=1]{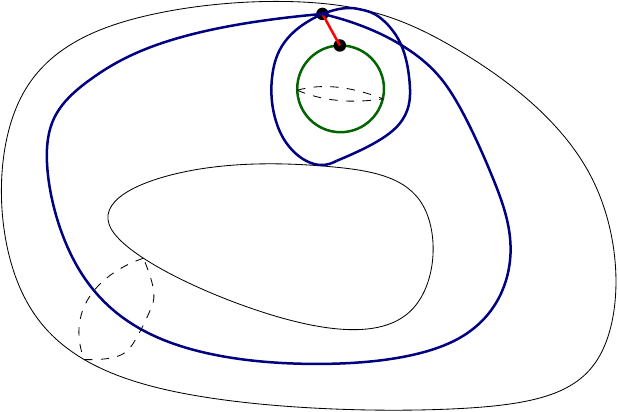}
            \caption{Cell decomposition of $N'$, with two vertices and four 1-cells. Cells of higher dimension are not displayed. The 2-cells correspond to: the sphere, the annulus which is visible around the sphere, and two more cells on the torus surface.}\label{fig:cell}
        \end{figure}
        We are going to define a homotopy $H \colon N' \times I \to \pi(N')$ between $\pi \circ f|_{N'}$ and $f^\vee \circ \pi |_{N'}$, where the domain has the obvious CW-complex structure induced by the one of $N'$, namely to each $i$-cell~$\sigma$ of $N'$ there is a corresponding $(i+1)$-cell~$\sigma \times I$ in $N' \times I$. 
        Since the $0$-skeleton of $N'$ entirely lies in the boundary of $N'$ we can define the homotopy to be stationary at $(N')^{(0)} \times I$. Similarly, for each $1$-cell~$e$ of $N'$ lying in $\partial N$ we define $H$ to be stationary at $e \times I$.
        For the unique $1$-cell of $N'$ that does not lie in $\partial N'$ (that is the red one in Figure~\ref{fig:cell}) we observe that the maps $u_0$ and $u_1$, restricted to it, give homotopic curves (relative to endpoints): They are both curves that, first, head to the wedge-sum basepoint $x_0$, and then travel through $\varphi\circ\alpha$.
        Thus, we define $H$ over $e \times I$ to be such a homotopy.
        The problem of extending $H$ to the remaining cells corresponds to being able to extend maps from the $m$-sphere ($m \geq 2$) $S^m \to \pi(N')$ to maps from the $(m+1)$-disc $D^{m+1} \to \pi(N')$. Since we have just proved that $\pi(N')$ is aspherical, this filling problem can be always successfully solved. This concludes the construction of the desired homotopy, whence the proof.
    \end{proof}
    
    Let $M^\vee_{f^\vee}$ be the mapping torus of $f^\vee$. By Lemma~\ref{lemma:first:homotopy:square:commute}, we know that Diagram~\eqref{eq:diagram_f} commutes up to homotopy.
    Hence, we have a continuous map~$\pi_f \colon M_f \to M^\vee_{f^\vee}$ (as explained in Subsection~\ref{subsection:digram:induces:map}). Moreover, since $\pi_* \colon \pi_1(M, w) \to \pi_1(M^\vee, x_0)$ is an isomorphism of fundamental groups, the induced map $\pi_f$ between the mapping tori is also a $\pi_1$-isomorphism (Proposition~\ref{prop:iso_fundgroups}).
    Finally, by Remark~\ref{rem:lifted:fund:classes} the facts that $H_3(f)([M]) = [M]$, $H_3(\pi)([M]) = [M^\vee]$ and, consequently, $H_3(f^\vee)([M^\vee]) = [M^\vee]$, shows that the fundamental classes of $M$ and $M^\vee$ can be both lifted to homology classes of the corresponding mapping tori: $[M_f] \in H_4(M_f;\mathbb{Z})$ and $[M^\vee_{f^\vee}] \in H_4(M^\vee_{f^\vee};\mathbb{Z})$. 
    The former, being $M_f$ an oriented closed connected manifold, is the usual fundamental class of $M_f$. The commutativity of Diagram~\eqref{eq:diagram_f} and Remark~\ref{rem:fundamental:class:goes:to:fund:class} show that $H_4(\pi_f)([M_f]) = [M_{f^\vee}^\vee]$. 

    The previous discussion together with \Cref{prop:properties:omega}.(\emph{5}) shows that a sufficient condition for the vanishing of $\minent(M_f)$ is the vanishing of $\Omega([M^\vee_{f^\vee}])$. This will be achieved in the next step which will conclude the proof of Theorem~\ref{thm:minent_reducible}.
    

\medskip

\textbf{(Step 2)\quad}
We construct a map $f^\sqcup\colon M^\sqcup\to M^\sqcup$ such that the diagram 
\begin{equation}\label{eq:diagram_fsq}\begin{tikzcd}
        M^\sqcup \arrow["f^\sqcup", r] \arrow["\iota", d] & M^\sqcup \arrow["\iota", d]\\
        M^\vee \arrow["f^\vee", r] & M^\vee
\end{tikzcd}\end{equation}
commutes up to homotopy.
Again, this is what we need, in virtue of ``functoriality" of mapping tori (see Subsection~\ref{subsection:digram:induces:map}) in order to show that $\iota \colon M^\sqcup \to M^\vee$ induces a map between the corresponding mapping tori.

Let us write $f^\vee = g \circ h$, where $h \colon M^\vee \to M^\vee$ is the composition of the self-maps of $M^\vee$ corresponding to diffeomorphisms of types \emph{1}, \emph{2} and \emph{3}, while $g \colon M^\vee \to M^\vee$ is the composition of the remaining ones, corresponding to sliding diffeomorphisms of $M$, as described in the previous step.
Both $h$ and $g$ are allowed to be the identity maps when there are no diffeomorphisms of the corresponding type in the factorization of $f$.

Associated to $h \colon M^\vee \to M^\vee$ there is a permutation $\sigma$ of the set $\{1, \cdots, k\}$ that witnesses the fact that the restriction of $h$ to each $M_i$ is an orientation-preserving (pointed at $x_0$) diffeomorphism 
onto the factor $M_{\sigma(i)}$ (\Cref{rem:elementary_collar}). For convenience, we denote such restrictions by $h_i = h |_{M_i} \colon M_i \to M_{\sigma(i)}$. 

We consider the self-diffeomorphism~$f^\sqcup \colon M^\sqcup \to M^\sqcup$ that on each component $M_i$ agrees with the map $h_i$. 

\begin{lem}\label{lem:homotopy_sq}
    Diagram \eqref{eq:diagram_fsq} commutes up to homotopy.
\end{lem}
\begin{proof}
    Since $M^\sqcup$ is disconnected, we have to show that for each $i \in \{1, \cdots, k\}$ the maps $f^\vee |_{M_i}$ and $\iota |_{M_{\sigma(i)}} \circ h_i$ are homotopic as maps from $M_i$ to $M^\vee$.
    If we write $f^\vee |_{M_i} = g |_{M_{\sigma(i)}} \circ h_i$, it is readily seen that it is sufficient to show that $g |_{M_{\sigma(i)}}$ is homotopic to the inclusion $\iota |_{M_{\sigma(i)}}$ as a map from $M_{\sigma(i)}$ to $M^\vee$.

    Denote by $g_{\sigma(i)}\colon M_{\sigma(i)} \to M^\vee$ the restriction of $g$ to $M_{\sigma(i)}$.
    We are going to give an explicit description of $g_{\sigma(i)}$ that will show that it is homotopic to the inclusion.

    We recall that in the situation of Setup~\ref{setup:reducible:manifolds}, we have fixed some collars of $\partial W$. We also recall that given a collar $C_{\sigma(i)} \subset M_{\sigma(i)}' \subset M$, the quotient $\pi(C_{\sigma(i)}) \subset M_{\sigma(i)} \subset M^\vee$ is diffeomorphic to a $3$-dimensional ball that is the cone over the boundary sphere (the component of $C_{\sigma(i)}$ not lying in $\partial W$) with vertex $x_0$.
    
    We are going to prove by induction that $g_{\sigma(i)}$ has the following explicit description:
    There exist an integer~$p \geq 0$ and a curve~$\beta\colon[1/2^p,1]\to M^\vee$ such that:
    \begin{itemize}
        \item The endpoints $\beta(1/2^p)$ and $\beta(1)$ coincide with the wedge-sum basepoint~$x_0$ of~$M^\vee$;
        \item $g_{\sigma(i)}$ is the identity outside $\pi(C_{\sigma(i)})$;
        \item $g_{\sigma(i)}(\pi(x,t)) = \pi(x,2^pt)$ for $(x,t) \in S^2\times[0,1/2^p]\subset C_{\sigma(i)}$;
        \item $g_{\sigma(i)}(\pi(x,t)) = \beta(t)$ for $(x,t) \in S^2\times[1/2^p,1]\subset C_{\sigma(i)}$.
    \end{itemize}
    The claim is proved by induction on the number of slide diffeomorphisms that appear in the factorization of $f$.
    If this number is $0$, then the claim holds with $p=0$ and $\beta(1) = x_0$.

    Inductively, suppose that $g_{\sigma(i)}\colon M_i \to M^\vee$ satisfies the claim, with $p$ and $\beta$ as above.
    Then, consider a composition $s \circ g_{\sigma(i)}$, where $s$ is a self-map of $M^\vee$ corresponding to a slide diffeomorphism of $M$ as described in the construction of $f^\vee$.
    There are two cases:
    \begin{itemize}
        \item If the slide diffeomorphism is a slide of a component $M_j \neq M_{\sigma(i)}$, then $s$ agrees with the identity on $M_{\sigma(i)}$, so that $s \circ g_{\sigma(i)}$ satisfies the claim with the same integer $p$ and curve $s\circ\beta$ obtained by composition of continuous maps (since $s(x_0) = x_0$);
        \item Otherwise, recall that $s$ is defined using a curve $\bar\alpha:[0,1]\to M^\vee$, with $\bar\alpha(0)=\bar\alpha(1)=x_0$, so that $s(\pi(x,t)) = \bar\alpha(2t-1)$ for $(x,t) \in S^2\times [1/2,1] \subset C_{\sigma(i)}$, while $s$ is the identity outside $\pi(C_{\sigma(i)})$ and $s(\pi(x,t)) =\pi(x,2t)$ for $(x,t) \in S^2\times[0,1/2]\subset C_{\sigma(i)}$.
        Then, $s \circ g_{\sigma(i)}$ satisfies the claim with integer parameter $p+1$ and curve $[1/2^{p+1},1] \to M^\vee$ sending
        \[\begin{cases}
            t \mapsto \bar\alpha(2^{p+1}t-1) &t \in [1/2^{p+1},1/2^p]\\
            t \mapsto s(\beta(t)) &t \in [1/2^p,1].
        \end{cases}\]
    \end{itemize}
    Now that we have an explicit description of $g_{\sigma(i)}$ we are ready to provide the desired homotopy $\theta \colon M_{\sigma(i)}\times [0,1] \to M^\vee$ between $g_{\sigma(i)}$ and the inclusion $\iota |_{M_{\sigma(i)}}$. We define $\theta$ as follows:
    \begin{itemize}
        \item $\theta(p,z) = p$ for every $p \not\in \pi(C_{\sigma(i)})$ and $z \in [0,1]$;
        \item $\theta(\pi(x,t),z) = g_{\sigma(i)}(\pi(x,(1-z+z/2^p)t))$ for $(x,t) \in S^2\times[0,1] = C_{\sigma(i)}$.
    \end{itemize}
    This map is well defined because for $t = 1$ the value of the formula defining $\theta(\pi(x,t),z)$ does not depend on $x$:
    \[ \theta(\pi(x,1),z) = g_{\sigma(i)}\left(\pi\left(x,\left(1-z+\frac{z}{2^p}\right)\cdot 1\right)\right) = \beta\left(1-z+\frac{z}{2^p}\right). \]
    Moreover, $\theta$ is continuous because for all $\pi(x, 0) \in \pi(C_{\sigma(i)})$ and $z \in [0, 1]$ we have
    \[
    \theta(\pi(x,0),z) = g_{\sigma(i)}\Big(\pi(x, 0)\Big) = \pi(x, 0) = (x, 0) \in S^2 \times \{0\},
    \] 
    agreeing with the fact that outside $\pi(C_{\sigma(i)})$ all points are fixed.
    
    Finally, $\theta$ is a homotopy between $g_{\sigma(i)}$ and $\iota |_{M_{\sigma(i)}}$ because:
    \begin{itemize}
        \item For $z = 0$, we have $\theta(p,0) = g_{\sigma(i)}(p)$ for every $p \in M_{\sigma(i)}$;
        \item For $z = 1$, we have for every $(x, t) \in S^2\times[0,1] = C_{\sigma(i)}$ that
        \[
        \theta(\pi(x,t),1) = g_{\sigma(i)}\left(\pi\left(x,\frac{t}{2^p}\right)\right) = \pi(x,t),
        \]
    \end{itemize}
    concluding the proof.
\end{proof}

We are now ready to conclude the proof of Theorem~\ref{thm:minent_reducible}.
Let $R \subset \{1, \cdots, k\}$ be a subset containing a representative of each orbit of the action by the permutation $\sigma$ on $\{1, \cdots, k\}$ and let $|r|$ denote the cardinality of the orbit containing $r \in R$. Then, the mapping torus of $f^\sqcup \colon M^\sqcup \to M^\sqcup$ that we denote by $M^\sqcup_{f^\sqcup}$ is the disjoint union of the following mapping tori:
\[
\bigsqcup_{r \in R} (M_{r})_{\rho_r},
\]
where $\rho_r = h_{\sigma^{|r|-1}(r)} \circ \dots \circ h_{{\sigma(r)}}\circ h_r$ is an orientation-preserving self-\hspace{0pt}diffeomorphism of $M_r$ (Corollary~\ref{cor:MT:disconnected:general}).
Hence, each connected component is a mapping torus over a prime manifold with respect to an orientation-preserving self-diffeomorphism.

Since  $H_3(f^\sqcup)([M^\sqcup]) = [M^\sqcup]$, by Remark~\ref{rem:lifted:fund:classes} we have a lifted homology class $[M_{f^\sqcup}^\sqcup]$ on $M^\sqcup_{f^\sqcup}$, which in this case is just the sum of the fundamental classes of the connected components $(M_{r})_{\rho_r}$.
Hence, $[M_{f^\sqcup}^\sqcup] = \sum_{r \in R} [(M_{r})_{\rho_r}]$.

By Lemma~\ref{lem:homotopy_sq} there exists a map $\iota_{f^\sqcup}\colon M_{f^\sqcup}^\sqcup \to M_{f^\vee}^\vee$ and by Remark~\ref{rem:fundamental:class:goes:to:fund:class} we also know that
$H_4(\iota_{f^\sqcup})([M_{f^\sqcup}^\sqcup]) = [M_{f^\vee}^\vee]$.
Moreover, by linearity we have that
\[[M_{f^\vee}^\vee] = \sum_{r \in R} H_4(\iota_{f^\sqcup})([(M_{r})_{\rho_r}]).\] 

We are now ready to compute $\Omega([M_{f^\vee}^\vee])$:
\begin{align*}
    \Omega([M_{f^\vee}^\vee]) &= \Omega\left(\sum_{r \in R} H_4(\iota_{f^\sqcup})([(M_{r})_{\rho_r}])\right) \\
    &\leq \sum_{r \in R} \Omega\big(H_4(\iota_{f^\sqcup})([(M_{r})_{\rho_r}])\big) \\
    &\leq \sum_{r \in R} \Omega([(M_{r})_{\rho_r}]) \\
    &= \sum_{r \in R} \minent\big((M_{r})_{\rho_r}\big)^4 \\
    &=0,
\end{align*}
where we use, in order, items~(\emph{1}), (\emph{2}) and (\emph{3}) of Proposition~\ref{prop:properties:omega}, and the last equality is a consequence of Corollary~\ref{cor:minent:MT:prime}.
This concludes the proof of \Cref{thm:minent_reducible}.

\appendix

\section{Mapping tori}\label{sec:prelim:MT}

We collect in this appendix some topological and group-theoretic properties of mapping tori. Some of the results contained in the appendix about the homology and the fundamental group of mapping tori are classical, but we could not find suitable references. For this reason we prefer to state (and prove) all the results we need in the paper.

We begin with the well-known definition of the mapping torus of a continuous map.
\begin{defi}[mapping torus]
    Let $X$ be a topological space and let $f\colon X\to X$ be a continuous map.
    The \emph{mapping torus of $f$}, which we denote by $X_f$, is the quotient of $X\times[0,1]$ by the finest equivalence relation such that $(x,1) \sim (f(x),0)$ for every $x \in X$.
\end{defi}

We consider $X$ as a subspace of $X_f$ by identifying it with $X\times\{0\}$. Note that this is a closed embedding.

\begin{rem}[mapping tori over manifolds]\label{rem:MT:category:smooth:top}
    If $M$ is a manifold of dimension $n$ and $f\colon M\to M$ is a homeomorphism, then $M_f$ is a manifold of dimension $n+1$.
Moreover, if $M$ is orientable and $f$ is orientation-preserving, then $M_f$ is orientable~\cite[Proposition~42.9]{friedl}.
If $M$ is a smooth manifold and $f$ is a diffeomorphism, then $M_f$ naturally acquires the structure of a smooth manifold~\cite[Proposition~49.3]{friedl}.
\end{rem}

\begin{rem}[mapping tori and equivalences]\label{rem:isotopic:map:give:isotopic:mapping:tori}
   Let $X$ be a CW-complex and let $f, g \colon X \to X$ be two continuous self-maps of $X$. If $f$ and $g$ are homotopic then $X_f$ is homotopy equivalent to $X_g$~\cite[Proposition~0.18]{hatcher}.

   Moreover, if $M$ is an oriented compact smooth manifold with (possibly empty) boundary and $f, g \colon M \to M$ are two orientation-preserving self-homeomorphisms (resp. self-diffeomorphisms) of $M$ then $M_f$ is homeomorphic to $M_g$ (resp. diffeomorphic)~\cite[Lemma~49.4]{friedl}. The standard references only cover the closed case, however the case of compact manifolds with non-empty boundary is completely analogous: The main difference is that $M \times [0,1]$ is a manifold \emph{with corners} and so we have to use smooth collars for manifolds with boundary. More precisely, we first have to notice that the glueing sites $M\times \{0\}$ and $M\times \{1\}$ are \emph{smooth boundaries} of $M \times I$~\cite[p. 1209]{friedl}. Then the mapping torus $M_f$ also admits a structure of manifold with corners, unique up to diffeomorphisms isotopic to the identity~\cite[Proposition~54.13]{friedl}. Note that since the corners $\partial M\times \{0,1\}$ of $M\times [0,1]$ are completely contained in the glueing locus, the mapping torus  $M_f$ has empty corner \cite[Proposition~54.13(2b)]{friedl} or, in other words, it is just a manifold with boundary. Moreover, also in this setting we have that the diffeomorphism type of $M_f$ depends only on the isotopy type of $f$. Indeed, the proof for the closed case~\cite[Lemma~49.4]{friedl} can be adapted to the setting of manifolds with corners substituting the usual collar of the boundary with the collar of the smooth boundary used when we glue manifolds with corners.
\end{rem}

\subsection{Maps between mapping tori}\label{subsection:digram:induces:map}
Suppose that $X$ and $Y$ are two topological spaces, and that we have continuous maps $f\colon X\to X$ and $g\colon Y\to Y$ so that we can consider the respective mapping tori $X_f$ and $Y_g$.
Suppose now that we have a map $\varphi\colon X \to Y$ such that the diagram
    \begin{equation}\label{eq:diagramxy}\begin{tikzcd}
        X \arrow["f", r] \arrow["\varphi", d] & X \arrow["\varphi", d]\\
        Y \arrow["g", r] & Y
    \end{tikzcd}\end{equation}
commutes up to homotopy, i.e., $\varphi\circ f$ and $g\circ\varphi$ are homotopic maps; let $H\colon X\times[0,1] \to Y$ be a homotopy, with $H(x,0) = g(\varphi(x))$ and $H(x,1) = \varphi(f(x))$.
Then, we have a continuous map from $X_f$ to $Y_g$, obtained by considering the map $\hat H\colon X\times[0,1] \to Y\times[0,1]$ defined by
\[\hat H(x,t) = (H(x,t),t)\]
and passing to the quotients.
This map depends on the chosen homotopy, and we denote it in this section with $\varphi_H\colon X_f \to Y_g$.
However, in other parts of the paper the dependence from the homotopy is not relevant for our discussion.

\begin{rem}[induced map extends the composition]
    The map $\varphi_H$ extends the map $g\circ\varphi\colon X\to Y$ (not the map $\varphi$), where as we already mentioned we identify $X$ and $Y$ with the leaves at height $0$ of the respective mapping tori.
\end{rem}

\begin{rem}[an alternative description when the diagram strictly commutes]
    If Diagram~\eqref{eq:diagramxy} strictly commutes, one could also define -- perhaps more naturally -- a map $\varphi'\colon X_f \to Y_g$ by setting
    \[\varphi'(x,t) = (\varphi(x),t)\]
    and passing to the quotients.
    This does not give the same result as the above more general procedure: If $H$ is the stationary homotopy between the equal maps $\varphi\circ f$ and $g\circ\varphi$, then $\varphi_H = g'\circ\varphi'$,
    where $g'\colon Y_g\to Y_g$ acts as $g$ separately on every ``leaf'' of the mapping torus.
    However, it is readily seen that $g'$ is homotopic to the identity, and therefore $\varphi_H$ is homotopic to $\varphi'$. 
\end{rem}

\subsection{Lifting homology classes}\label{sec:fund_classes}
If $M$ is an oriented closed manifold and $f$ is an orientation-preserving homeomorphism of $M$, then the mapping torus $M_f$ is again an oriented closed manifold, and as such it has a fundamental class in its top-degree homology.

For most of the applications, it is convenient to consider spaces $X$ which are not necessarily manifolds, but finite dimensional CW-complexes, on which we will still have some ``preferred'' classes in homology that we wish to ``lift'' to the mapping tori.
This can be done using the following proposition (we omit the module of coefficients; we need only integer coefficients in other sections but any Abelian group can be used here).

\begin{prop}[functorial LES of mapping tori~{\cite[p.\ 5]{bucherneofytidis}}]\label{prop:mapping_es}
    Let $X$ be a topological space and $f\colon X\to X$ be a continuous map. Let $X_f$ be the mapping torus of $f$.
    Then, there is a long exact sequence
    \begin{equation*}
        \begin{tikzcd}
            \cdots \ar[r] & H_k(X) \ar[r] & H_k(X) \ar[r] & H_k(X_f) \ar[r] & H_{k-1}(X)\ar[r] & \cdots
        \end{tikzcd}
    \end{equation*}
    where the leftmost map is $H_k(f)-\mathrm{id}_{H_k(X)}$.
    This sequence is natural in the sense that, if $g\colon Y \to Y$ is a continuous map and $\varphi\colon X\to Y$ makes Diagram \eqref{eq:diagramxy} commute up to homotopy, then the diagrams
    \begin{equation}\label{diagram:two:sequences}
        \begin{tikzcd}
            H_k(X) \ar[r]\ar[d] & H_k(X) \ar[r]\ar[d] & H_k(X_f) \ar[r]\ar[d] & H_{k-1}(X)\ar[d]\\
            H_k(Y) \ar[r] & H_k(Y) \ar[r] & H_k(Y_g) \ar[r] & H_{k-1}(Y)
        \end{tikzcd}
    \end{equation}
    commute, where $Y_g$ is the mapping torus of $g$ and the vertical maps are induced by $g\circ\varphi\colon X\to Y$ and $\varphi_H\colon X_f\to Y_g$, respectively; the latter being the map induced by any homotopy $H$ between $g\circ\varphi$ and $\varphi\circ f$.
\end{prop}
\begin{proof}
    The construction of the long exact sequence is described in the book of Hatcher ~\cite[page 151]{hatcher}.
    The statement about commutativity follows from the construction.
\end{proof}

Suppose now that $H_{n+1}(X) = 0$, e.g., this holds when $X$ is a manifold or a CW-complex of dimension $n$, and that there exists a homology class $\alpha \in H_n(X)$ such that $H_n(f)(\alpha) = \alpha$.
Then, the exact sequence in \Cref{prop:mapping_es} implies the existence of a unique class $\alpha_f \in H_{n+1}(X_f)$ which is sent to $\alpha \in H_n(X)$; this is how we ``lift'' homology classes on mapping tori.

\begin{rem}[fundamental classes lift to fundamental classes]\label{rem:lifted:fund:classes}
    In the case where $f$ is an orientation-preserving self-\hspace{0pt}homeomorphism of an oriented (not necessarily connected) manifold $M$, and $\alpha = [M]$ is the fundamental class, then the class $\alpha_f \in H_{n+1}(M_f)$ obtained as above is the fundamental class of $M_f$.
\end{rem}

\begin{rem}[functoriality of the lifted classes]\label{rem:fundamental:class:goes:to:fund:class}
    If $X$ and $Y$ are topological spaces with vanishing homology in degree $n+1$, and we have a diagram of continuous maps \eqref{eq:diagramxy} commuting up to homotopy, and we also have classes $\alpha \in H_n(X)$ and $\beta \in H_n(Y)$ such that $H_n(f)(\alpha) = \alpha$ and $H_n(g)(\beta) = \beta$, then we have $H_{n+1}(\varphi_H)(\alpha_f) = \beta_f$, because of the commutativity statement in \Cref{prop:mapping_es}.
\end{rem}

\subsection{Fundamental groups of mapping tori}
We believe that the facts that we now discuss about the fundamental group of mapping tori should be known to experts, but we include them here since we have not been able to find appropriate references in the literature.

Suppose that $X$ is a path-connected topological space, and fix a basepoint $x_0 \in X$.
Let $f\colon X\to X$ be a continuous map.
Upon chosing a path from $f(x_0)$ to $x_0$, which allows to identify $\pi_1(X,x_0) \cong \pi_1(X,f(x_0))$, we have an induced homomorphism $f_* \colon \pi_1(X,x_0)\to \pi_1(X,x_0)$.

We denote the chosen path with $u_X$; hence, if an element of $\pi_1(X,x_0)$ is represented by a path $\gamma$, we have that $f_*[\gamma] = [\overline{u_X}\star (f\circ\gamma) \star u_X]$, where the overline denotes the reverse path and $\star$ is the usual path concatenation.

Let $X_f$ be the mapping torus of $f$.
Then, $X_f$ is path-connected and its fundamental group $\pi_1(X_f,x_0)$ is isomorphic to the group
\begin{equation}\label{eq:fund_group}
\pi_1(X_f,x_0) \cong (\pi_1(X,x_0) * \langle s \rangle) / R_f,
\end{equation}
where $\langle s \rangle$ is an infinite cyclic group generated by $s$ and $R_f$ is the normal subgroup of the free product $\pi_1(X,x_0) * \langle s \rangle$ normally generated by the relations
\[ \{s^{-1}[\gamma] s = f_*[\gamma] \mid [\gamma] \in \pi_1(X,x_0)\}. \]
Through this explicit isomorphism describing $\pi_1(X_f, x_0)$, the element $s$ is represented by the path $v_X\star u_X$, where $v_X$ is the ``vertical'' path from $x_0 \in X = X\times\{0\} \subset X_f$ to $(x_0,1) = (f(x_0),0)$.
The relations in $R_f$ express the fact that, for every $[\gamma] \in \pi_1(X,x_0)$, the path $\gamma$ is homotopic (in $X_f$) to $v_X\star(f\circ\gamma)\star \overline{v_X}$.

\begin{rem}[fundamental group of mapping tori of $\pi_1$-automorphisms]\label{rem:mapping:torus:group:semidirect}
    If $f\colon X\to X$ induces an automorphism of the fundamental group, then the fundamental group of the resulting mapping torus is the semidirect product $\pi_1(X,x_0) \rtimes_{f_*} \Z$.
    In the general case it is a semidirect product of a \emph{quotient} of $\pi_1(X,x_0)$.
\end{rem}

Suppose now that we have another topological space $Y$ and maps $g\colon Y\to Y$ and $\varphi\colon X\to Y$ so that Diagram \eqref{eq:diagramxy} commutes up to homotopy~$H$, and consequently we also have a continuous map $\varphi_H\colon X_f \to Y_g$.
Let $y_0 = g(\varphi(x_0)) = \varphi_H(x_0)$.
We wish to describe the induced homomorphism
\[(\varphi_H)_*\colon \pi_1(X_f,x_0) \to \pi_1(Y_g,y_0)\]
as a map
\[ (\pi_1(X,x_0) * \langle s \rangle) / R_f \to (\pi_1(Y,y_0) * \langle s \rangle) / R_g.\]

The precise description depends on the chosen paths, $u_X$ and $u_Y$, joining $f(x_0)$ to $x_0$ in $X$ and $g(y_0)$ to $y_0$ in $Y$, respectively, which are needed in order to describe the fundamental groups as in~\eqref{eq:fund_group}.
Let us describe a convenient way to make these choices.
\begin{itemize}
    \item First, choose an arbitrary path $u_X$ from $f(x_0)$ to $x_0$ in $X$;
    \item Let $z$ be the path in $Y$ defined by $z(t)=H(x_0,t)$, which starts at $y_0$ and ends in $\varphi(f(x_0))$.
    Then, define $w = z \star (\varphi\circ u_X)$, which starts at $y_0$ and ends in $\varphi(x_0)$.
    Finally, take $u_Y = g\circ w$: This is a path in $Y$ from $g(y_0)$ to $g(\varphi(x_0)) = y_0$.
\end{itemize}

Now, the map $(\varphi_H)_*$ has the following description:
It is induced (passing to the quotients) by the map
\[ \pi_1(X,x_0) * \langle s \rangle \to \pi_1(Y,y_0) * \langle s \rangle\]
which sends $s$ to $s$ and sends any $[\gamma] \in \pi_1(X,x_0)$ to $(g\circ\varphi)_*[\gamma] \in \pi_1(Y,y_0)$.
The fact that $s\mapsto s$ is due to the fact that the path $\varphi_H\circ v_X$ is homotopic (relative to the endpoints) to $v_Y \star (g\circ z)$.

Consider now the map between fundamental groups induced by $\varphi$.
Strictly speaking, $\varphi_*$ is a homomorphism from $\pi_1(X,x_0)$ to $\pi_1(Y,\varphi(x_0))$.
We identify the latter with $\pi_1(Y,y_0)$ using the path $w$ defined above, so that we consider $\varphi_*$ as a map $\varphi_*\colon \pi_1(X,x_0)\to \pi_1(Y,y_0)$ sending any $[\gamma]\in\pi_1(X,x_0)$ to $[w\star (\varphi\circ\gamma) \star\overline{w}]$.
With this choice, we have by construction that $(g\circ\varphi)_* = g_* \circ \varphi_*$.
Moreover,
\begin{align*}
    \varphi_*f_*[\gamma] &= \varphi_*[\overline{u_X}\star (f\circ\gamma)\star u_X] \\
    &= [w\star (\varphi\circ \overline{u_X})\star(\varphi\circ f\circ\gamma)\star(\varphi\circ u_X)\star\overline{w}]\\ 
    &= [z \star (\varphi\circ f\circ\gamma)\star\overline{z}] \\
    &= [g\circ\varphi\circ\gamma],
\end{align*}
so that we also have $\varphi_*\circ f_* = g_*\circ\varphi_*$ (the last equality follows by considering the map $(t_1,t_2)\mapsto H(\gamma(t_1),t_2)$ defined on $I\times I$ and inspecting its restriction to the boundary of $I\times I$). 

\begin{prop}[$\varphi$ $\pi_1$-iso $\Rightarrow \varphi_H$ $\pi_1$-iso]\label{prop:iso_fundgroups}
    Let $f,g$ and $\varphi$ be continuous maps such that Diagram~\eqref{eq:diagramxy} commutes up to a homotopy~$H$. Let $\varphi_H\colon X_f \to Y_g$ be the induced map between the mapping tori.
    If $\varphi_*$ is an isomorphism, then $(\varphi_H)_*$ is an isomorphism.
\end{prop}

\begin{proof}
    Let $c\colon \pi_1(Y_g,y_0) \to \pi_1(Y_g,y_0)$ be the conjugation $c(h) = sh s^{-1}$.
    Since $c$ is an automorphism, to conclude it suffices to show that the composition $c \circ (\varphi_H)_*$ is an isomorphism.
    This composition is induced by the map
    \[ \pi_1(X,x_0) * \langle s \rangle \to \pi_1(Y,y_0) * \langle s \rangle\]
    sending $s$ to $s$ and any $[\gamma] \in \pi_1(X,x_0)$ to $s\cdot (g\circ\varphi)_*[\gamma]\cdot s^{-1} \in \pi_1(Y,y_0)*\langle s\rangle$.
    Using the relations in $R_g$, we have that the computation
    \begin{align*}
        s\cdot (g\circ\varphi)_*[\gamma]\cdot s^{-1} &= s \cdot g_*(\varphi_*[\gamma]) \cdot s^{-1} \\
        &= s \cdot s^{-1} \cdot \varphi_*[\gamma] \cdot s \cdot s^{-1}\\
        &= \varphi_*[\gamma]
    \end{align*}
    holds in the quotient $\pi_1(Y,y_0)*\langle s\rangle / R_g$.
Then the map $c\circ (\varphi_H)_*$ is also induced by the map
\[\psi\colon\pi_1(X,x_0) * \langle s \rangle \to \pi_1(Y,y_0) * \langle s \rangle \]
such that $\psi(s)=s$ and $\psi([\gamma])=\varphi_*[\gamma]$. It is easy to check that since $\varphi_*$ is an isomorphism, also $\psi$ is an isomorphism. Moreover, since $\varphi_*\circ f_* = g_*\circ\varphi_*$,  for every $[\gamma]\in \pi_1(X,x_0)$ we have that
\begin{align*}
    \psi(s^{-1}\cdot [\gamma] \cdot s\cdot (f_*[\gamma])^{-1})&=s^{-1}\cdot \varphi_*[\gamma]\cdot s \cdot (g_*(\varphi_*[\gamma]))^{-1}.
\end{align*}
Since $\psi$ is a isomorphism, this readily implies that $\psi(R_f)= R_g$, concluding the proof. 
    %
\end{proof}

\subsection{Mapping tori of composition of maps}

In this section we give an explicit description of mapping tori of self-homeomorphisms of disconnected topological spaces as disjoint unions of mapping tori over connected spaces. 

\begin{setup}\label{setup:composition:MT}
    We consider the following setup:
    \begin{itemize}
        \item Let $X$ be a connected topological space and let $n \geq 2$ be an integer;

        \item Let $X_1, \cdots, X_n$ be connected topological spaces homeomorphic to $X$. For each $i \in \{1, \cdots, n\}$ we fix a homeomorphism $\varphi_i \colon X \to X_i$. Moreover, for every $x\in X$, denote by $x_i\in X_i$ the image of $x$ via $\varphi_i$;

        \item For every $i\in \{1,\cdots,n\}$ we fix a homeomorphism $f_i\colon X_i\to X_{i+1}$ with the convention that $X_{n+1}=X_1$.
    \end{itemize}
\end{setup}
In the situation of Setup~\ref{setup:composition:MT} we consider the finest equivalence relation~$\sim$ on $\sqcup_{i=1}^n (X_i \times [0,1])$ such that
\[
\forall x \in X, \ \forall i\in \{1,\cdots, n\} \quad (x_i,1)\sim(f_i(x_i),0)
\]
and we denote by $MT(X;f_1,\cdots, f_n)$ the quotient space $(\sqcup_{i=1}^n (X_i \times [0,1]))/ \sim$.
By definition this is exactly the mapping torus of the map $f\colon \sqcup_{i=1}^nX_i\to \sqcup_{i=1}^nX_i$ such that $f|_{X_i}=f_i$.

\begin{prop}[mapping torus of composition of maps]\label{prop:MT:disconnected:cyclic}
    In the situation of Setup~\ref{setup:composition:MT}, denote by $g\colon X_1\to X_1$ the composition $g=f_{n}\circ\dots\circ f_2\circ f_1$.
    Then, the following spaces are homeomorphic:
    \[
    MT(X;f_1,\cdots,f_n)\cong (X_1)_g.
    \]
\end{prop}
\begin{proof}
Let $\pi\colon X_1 \times [0,1]\to (X_1)_g$ be the quotient map and consider the continuous map $\psi\colon \sqcup_{i=1}^n (X_i \times [0,1])\to (X_1)_g$ defined by the formula
\[
h(x_i,t)=
\begin{cases}
\pi(x_i,\frac{t}{n}) & \text{if }i=1\\
   \pi\left((f_{i-1}\circ\dots\circ f_1)^{-1}(x_i), \frac{t+i-1}{n}\right) & \text{if }i\neq 1.
\end{cases}
\]
We have that for every $i \in\{2,\cdots, n\}$
\[ h(x_i,0)= \pi\left(f_1^{-1}\circ\dots\circ f_{i-1}^{-1}(x_i), \frac{i-1}{n}\right) =h(f_{i-1} ^{-1}(x_{i}),1) \]
and for $i=1$ we have
\[
h(x_1,0)=\pi(x_1,0)=\pi(f_1^{-1}\circ\dots\circ f_n^{-1}(x_1),1)=h(f_n^{-1}(x_1),1).
\]
This shows that $h$ is constant on each equivalence class, hence it induces a continuous map~$h'\colon MT(X;f_1,\cdots,f_n)\to (X_1)g$.

We are now going to provide a continuous bi-inverse of $h'$. First, denote by $q\colon\sqcup_{i=1}^n (X_i \times [0,1])\to MT(X;f_1,\cdots,f_n)$ the quotient map and consider the map $s\colon X_1 \times [0,1]\to  MT(X;f_1,\cdots,f_n)$ defined by the formula
\[
s(x_1,t)=\begin{cases}
    q(x_1,nt) &\text{if } t\in[0,\frac{1}{n}]\\
    q(f_{i}\circ\dots\circ f_1(x_1), nt-i) &\text{if }t\in [\frac{i}{n},\frac{i+1}{n}], \ i\in\{1,\cdots, n-1\}.
\end{cases}
\]
This is well defined since by construction for every $i \in \{1, \cdots, n-1\}$ we have
\[ q(f_{i-1}\circ \dots \circ f_1(x_1), 1))= q(f_i\circ \dots \circ f_1(x_1), 0). \]
Moreover, $s$ induces a continuous map $s'\colon (X_1)_g\to MT(X;f_1,\cdots,f_n)$ because of the following computation:
\begin{align*}
    s(x_1,1)&= q(f_{n-1}\circ\dots\circ f_1(x_1), 1)
    \\&= q(f_{n}\circ\dots\circ f_1(x_1), 0)\\&=q(g(x_1),0)\\&=s(g(x_1),0).
\end{align*}
It is immediate to check that $s'$ provides a bi-inverse of $h'$, whence they are both homeomorphisms and $(X_1)_g\cong MT(X;f_1,\cdots,f_n)$.
\end{proof}
\begin{rem}[the equivalence also works in the smooth category]\label{rem:MT:disconnected:smooth}
    When $X$ is a smooth manifold and the maps $f_i$ are diffeomorphisms, the maps $h'$ and $s'$ defined in the proof above are smooth. Hence we have that $MT(X;f_1,\cdots,f_n)$ and $(X_1)_g$ are diffeomorphic as smooth manifolds with respect to the induced smooth structure (Remark~\ref{rem:MT:category:smooth:top}). 
\end{rem}
Proposition~\ref{prop:MT:disconnected:cyclic} allows us to give an explicit description of mapping tori over disconnected spaces under the assumption that the number of connected components is finite.
\begin{setup}\label{setup:MT:composition:disconnected}
  We consider the following situation:
  \begin{itemize}
      \item Let $m \geq 1$ be an integer and suppose that $X = N_1 \sqcup \cdots \sqcup N_m$ is the decomposition of $X$ into its connected components; 

      \item Let $f \colon X \to X$ be a self-homeomorphism of $X$;

      \item  Up to reodering $\{N_1, \cdots, N_m\}$ we may assume that $\{N_1, \cdots, N_k\} \subset \{N_1, \cdots, N_m\}$ is a set of representative of these $f$-orbits. Since $m$ is assumed to be finite, $f$ cyclically permutes the connected components of each orbit;

      \item By construction for each $j \in \{1, \cdots, k\}$ there exists a minimal integer $n_j \geq 1$ such that $f^{n_j}(N_j)=N_j$. We denote by $X_j = \sqcup_{i = 1}^{n_j} f^i(N_j)$ and we have that $X = \sqcup_{j = 1}^k X_j$ is a decomposition of $X$ in $f$-invariant disjoint sets;

      \item For all $j \in \{1, \cdots, k\}$ and $i \in \{1, \cdots, n_j\}$ we denote by $g_j \colon N_j \to N_j$ the restriction of $f^{n_j}$ to $N_j$, by $f_j \colon X_j \to X_j$ the restriction of $f$ to $X_j$ and by $f_{i,j} \colon f^{i-1}(N_j) \to f^{i}(N_j)$ the restriction of $f$ to $f^{i-1}(N_j)$ (here we set $f^0$ to be the identity).
  \end{itemize}
\end{setup}
We have the following corollary:
\begin{cor}\label{cor:MT:disconnected:general}
    In the situation of Setup~\ref{setup:MT:composition:disconnected} the mapping torus~$X_f$ is homeomorphic to the disjoint union of mapping tori of $g_j$:
    \[
    X_f \cong \bigsqcup_{j=1}^k(N_j)_{g_j}.
    \]
    In particular, if $X$ is a compact smooth manifold and $f$ is a self-\hspace{0pt}diffeomorphism,  the analogue statement holds in the category of smooth manifolds.
\end{cor}
 \begin{proof}
    By definition for every $j \in \{1, \cdots, k\}$ we have that the mapping torus~$(X_j)_{f_j}$ is homeomorphic to the space~$MT(N_j;f_{1,j},\cdots, f_{n_j, j})$. Moreover, by Proposition \ref{prop:MT:disconnected:cyclic} the latter space is homeomorphic to the mapping torus~$(N_j)_{g_j}$. Since by construction the space $X_j$ is closed, open in $X$ and $f$-invariant, the resulting mapping torus $(X_j)_{f_j}$ is a connected component of $X_f$. Thus, we have proved that
    \[
    X_f \cong \bigsqcup_{j = 1}^k (N_j)_{g_j}.
    \]
     In the case where $X$ is a compact smooth manifolds and $f$ in a self-diffeomor\-phism, the analogous statement in the smooth category follows from Remark~\ref{rem:MT:disconnected:smooth}.
 \end{proof}

\subsection*{Declarations} There is no conflict of interest to declare. There is no relevant data related to the paper.

\bibliographystyle{alpha}
\bibliography{svbib}

\newcommand{\etalchar}[1]{$^{#1}$}
\begin{thebibliography}{GLGAH13}

\bibitem[AFW15]{aschenbrenner20123}
M.~Aschenbrenner, S.~Friedl, and H.~Wilton.
\newblock {\em 3-manifold groups}.
\newblock EMS Series of Lectures in Mathematics. European Mathematical Society
  (EMS), Z\"urich, 2015.

\bibitem[Bab93]{babenko1993asymptotic}
I.~K. Babenko.
\newblock Asymptotic invariants of smooth manifolds.
\newblock {\em Russian Acad. Sci. Izv. Math.}, 41(1):1--38, 1993.

\bibitem[BBF{\etalchar{+}}14]{bucher2014isometric}
M.~Bucher, M.~Burger, R.~Frigerio, A.~Iozzi, C.~Pagliantini, and M.~B.
  Pozzetti.
\newblock Isometric embeddings in bounded cohomology.
\newblock {\em J.\ Topol.\ Anal.}, 6(1):1--25, 2014.

\bibitem[BBM{\etalchar{+}}10]{francesi}
L.~Bessi{\`e}res, G.~Besson, S.~Maillot, M.~Boileau, and J.~Porti.
\newblock {\em Geometrisation of 3-manifolds}, volume~13 of {\em EMS Tracts in
  Mathematics}.
\newblock European Mathematical Society (EMS), {Z}{\"u}rich, 2010.

\bibitem[BC21]{bregman2021minimal}
C.~Bregman and M.~Clay.
\newblock Minimal volume entropy of free-by-cyclic groups and 2-dimensional
  right-angled artin groups.
\newblock {\em Math.\ Ann.}, 381(3):1253--1281, 2021.

\bibitem[BN20]{bucherneofytidis}
M.~Bucher and C.~Neofytidis.
\newblock The simplicial volume of mapping tori of $3$-manifolds.
\newblock {\em Math.\ Ann.}, 376(3-4):1429--1447, 2020.

\bibitem[BP92]{benedetti1992lectures}
Riccardo Benedetti and Carlo Petronio.
\newblock {\em Lectures on hyperbolic geometry}.
\newblock Springer Science \& Business Media, 1992.

\bibitem[Bru08]{brunnbauer2008homological}
M.~Brunnbauer.
\newblock Homological invariance for asymptotic invariants and systolic
  inequalities.
\newblock {\em Geom.\ Funct.\ Anal., GAFA}, 18(4):1087--1117, 2008.

\bibitem[BS24]{BS2023seminorm}
I.~Babenko and S.~Sabourau.
\newblock Volume entropy semi-norm and systolic volume semi-norm.
\newblock {\em J.\ Eur.\ Math.\ Soc. (JEMS)}, 26(11):4393--4439, 2024.

\bibitem[BS25]{BS2025fiber}
I.~Babenko and S.~Sabourau.
\newblock Minimal volume entropy and fiber growth.
\newblock {\em J.\ {\'E}c.\ polytech., Math.}, 12:481--521, 2025.

\bibitem[CG86]{cheegergromov}
Jeff Cheeger and Mikhael Gromov.
\newblock {Collapsing {R}iemannian manifolds while keeping their curvature
  bounded.~I}.
\newblock {\em J.~Differential~Geom.}, 23(3):309--346, 1986.

\bibitem[CLM22]{CLM}
Pietro Capovilla, Clara L{\"o}h, and Marco Moraschini.
\newblock Amenable category and complexity.
\newblock {\em Algebraic \& Geometric Topology}, 22(3):1417--1459, 2022.

\bibitem[CLN85]{camachoneto}
C.~Camacho and A.~Lins~Neto.
\newblock {\em Geometric theory of foliations}.
\newblock Birkh\"auser Boston, Inc., Boston, MA, 1985.
\newblock Translated from the Portuguese by Sue E. Goodman.

\bibitem[CZ06]{CaoZhu}
H.-D. Cao and X.~P. Zhu.
\newblock A complete proof of the {P}oincar{\'e} and geometrization conjectures
  -- application of the {H}amilton-{P}erelman theory of the {R}icci flow.
\newblock {\em Asian J. Math.}, 10:165--492, 2006.

\bibitem[Din71]{dinaburg1971connection}
E.~I. Dinaburg.
\newblock A connection between various entropy characterizations of dynamical
  systems.
\newblock {\em Izv. Akad. Nauk SSSR Ser. Mat}, 35(324-366):13, 1971.

\bibitem[Fri17]{frigerio2017bounded}
R.~Frigerio.
\newblock {\em Bounded cohomology of discrete groups}, volume 227 of {\em
  Mathematical Surveys and Monographs}.
\newblock American Mathematical Society, Providence, RI, 2017.

\bibitem[Fri24]{friedl}
S.~Friedl.
\newblock Topology.
\newblock
  \url{https://friedl.app.uni-regensburg.de/papers/1t-total-public-october-7-2024.pdf},
  2024.

\bibitem[GLGAH13]{GGH}
J.~G\'{o}mez-Larra{\~n}aga, F.~Gonz\'{a}lez-Acu{\~n}a, and W.~Heil.
\newblock Amenable category of three-manifolds.
\newblock {\em Algebr.\ Geom.\ Top.}, 13:905--925, 2013.

\bibitem[Gro81]{gromov1981groups}
M.~Gromov.
\newblock Groups of polynomial growth and expanding maps (with an appendix by
  jacques tits).
\newblock {\em Publ.\ Math.\ IH{\'E}S}, 53:53--78, 1981.

\bibitem[Gro82]{vbc}
M.~Gromov.
\newblock Volume and bounded cohomology.
\newblock {\em Inst. Hautes Études Sci. Publ. Math.}, 56:5--99 (1983), 1982.

\bibitem[Hat02]{hatcher}
A.~Hatcher.
\newblock {\em Algebraic topology}.
\newblock Cambridge University Press, Cambridge, 2002.

\bibitem[Hat07]{hatcher2007notes}
A.~Hatcher.
\newblock Notes on basic 3-manifold topology.
\newblock \url{https://pi.math.cornell.edu/~hatcher/3M/3Mfds.pdf}, 2007.

\bibitem[Hil02]{hillman2002four}
J.~A. Hillman.
\newblock {\em Four-manifolds, geometries and knots}, volume~5 of {\em Geometry
  \& Topology Monographs}.
\newblock Geometry \& Topology Publications, Coventry, 2002.

\bibitem[JWW01]{jiang1996homeomorphisms}
B.~Jiang, S.~Wang, and Y.-Q. Wu.
\newblock Homeomorphisms of 3-manifolds and the realization of nielsen number.
\newblock {\em Comm. Anal. Geom.}, 9(4):825--877, 2001.

\bibitem[KKM06]{kedra2006crossed}
J.~Kedra, D.~Kotschick, and S.~Morita.
\newblock Crossed flux homomorphisms and vanishing theorems for flux groups.
\newblock {\em Geom.\ Funct.\ Anal., GAFA}, 16(6):1246--1273, 2006.

\bibitem[KL08]{KleLott}
B.~Kleiner and J.~Lott.
\newblock Notes on {P}erelman's papers.
\newblock {\em Geom. Topol.}, 12:2587--2855--492, 2008.

\bibitem[KN07]{kapovich2007patterson}
I.~Kapovich and T.~Nagnibeda.
\newblock The {P}atterson--{S}ullivan embedding and minimal volume entropy for
  outer space.
\newblock {\em Geom.\ Funct.\ Anal., GAFA}, 17(4):1201--1236, 2007.

\bibitem[Kne29]{Kneser}
H.~Kneser.
\newblock Geschlossene {F}l{\"a}chen in dreidimensionalen {M}annigfaltigkeiten.
\newblock {\em Jahresbericht der Deutsche Mathematische Vereinigung},
  38:248--260, 1929.

\bibitem[Kot12]{kotschick2011entropies}
D.~Kotschick.
\newblock Entropies, volumes, and {E}instein metrics.
\newblock In {\em Global differential geometry}, volume~17 of {\em Springer
  Proc. Math.}, pages 39--54. Springer, Heidelberg, 2012.

\bibitem[Lim08]{Lim08}
S.~Lim.
\newblock Minimal volume entropy for graphs.
\newblock {\em Trans. Amer. Math. Soc.}, 360(10):5089--5100, 2008.

\bibitem[LM22]{lmfibration}
C.~L{\"o}h and M.~Moraschini.
\newblock Topological volumes of fibrations: a note on open covers.
\newblock {\em Proc.\ R.\ Soc.\ Edin.\ Sec.\ A: Math.}, 152(5):1340--1360,
  2022.

\bibitem[LMR22]{loehmoraschiniraptis}
C.~L\"oh, M.~Moraschini, and G.~Raptis.
\newblock On the simplicial volume and the euler characteristic of (aspherical)
  manifolds.
\newblock {\em Res.\ Math.\ Sci.}, 9(3):44, 2022.

\bibitem[LMS22]{LMS}
C.~L{\"o}h, M.~Moraschini, and R.~Sauer.
\newblock Amenable covers and integral foliated simplicial volume.
\newblock {\em New York J.\ Math.}, 28:1112--1136, 2022.

\bibitem[L{\"o}h20]{loehergodic}
C.~L{\"o}h.
\newblock {\em Ergodic Theoretic Methods in Group Homology: {A} Minicourse on
  {$L^2$}-{B}etti numbers in group theory}.
\newblock SpringerBriefs in Matehmatics. Springer, Cham, 2020.

\bibitem[Man79]{manning1979topological}
A.~Manning.
\newblock Topological entropy for geodesic flows.
\newblock {\em Ann.\ Math.}, 110(3):567--573, 1979.

\bibitem[Mar16]{martelli}
B.~Martelli.
\newblock {\em An introduction to geometric topology}.
\newblock 2016.

\bibitem[McC86]{McCull}
D.~McCullough.
\newblock Mappings of reducible 3-manifolds.
\newblock In {\em Proceedings of the semester in Geometric and Algebraic
  Topology, Warsaw, Banach Center}, pages 61--76, 1986.

\bibitem[McM15]{mcmullen2015entropy}
C.~T. McMullen.
\newblock Entropy and the clique polynomial.
\newblock {\em J. Topol.}, 8(1):184--212, 2015.

\bibitem[Mil62]{Milnor}
J.~Milnor.
\newblock A {U}nique {D}ecomposition {T}heorem for 3-{M}anifolds.
\newblock {\em American Journal of Mathematics}, 84(1):1--7, 1962.

\bibitem[Mil68]{milnor1968note}
J.~Milnor.
\newblock A note on curvature and fundamental group.
\newblock {\em J.\ Diff.\ Geom.}, 2(1):1--7, 1968.

\bibitem[MM03]{moerdijk-book}
I.~Moerdijk and J.~Mr\v{c}un.
\newblock {\em Introduction to foliations and {L}ie groupoids}, volume~91 of
  {\em Cambridge Studies in Advanced Mathematics}.
\newblock Cambridge University Press, Cambridge, 2003.

\bibitem[MN20]{mann2020dynamical}
K.~Mann and S.~Nariman.
\newblock Dynamical and cohomological obstructions to extending group actions.
\newblock {\em Math. Ann.}, 377(3):1313--1338, 2020.

\bibitem[MP99]{moerdijk+pronk}
I.~Moerdijk and D.~A. Pronk.
\newblock Simplicial cohomology of orbifolds.
\newblock {\em Indag. Math. (N.S.)}, 10(2):269--293, 1999.

\bibitem[MT07]{MorganTian}
J.~Morgan and G.~Tian.
\newblock {\em Ricci flow and the {P}oincar{\'e} conjecture}, volume~3 of {\em
  Clay Mathematics Monographs}.
\newblock American Mathematical Society, Providence, RI, Clay Mathematics
  Institute, Cambridge, MA, 2007.

\bibitem[Pera]{Per02}
G.~Perelman.
\newblock The entropy formula for the {R}icci flow and its geometric
  applications.
\newblock arXiv:math/0211159.

\bibitem[Perb]{Per03bis}
G.~Perelman.
\newblock Finite extinction time for the solutions to the {R}icci flow on
  certain three-manifolds.
\newblock arXiv:math/0307245.

\bibitem[Perc]{Per03}
G.~Perelman.
\newblock Ricci flow with surgery on three-manifolds.
\newblock arXiv:math/0303109.

\bibitem[Pet95]{pettet1995finitely}
M.~R. Pettet.
\newblock Finitely generated groups with virtually free automorphism groups.
\newblock {\em Proc.\ Edin.\ Math.\ Soc.}, 38(3):475--484, 1995.

\bibitem[Pie19]{pieroni2019minimal}
E.~Pieroni.
\newblock Minimal entropy of $3$-manifolds.
\newblock arXiv:1902.09190, 2019.

\bibitem[PP03]{paternainpetean}
G.~Paternain and J.~Petean.
\newblock {Minimal entropy and collapsing with curvature bounded from below}.
\newblock {\em Invent.~Math.}, 151:415--450, 2003.

\bibitem[Sam99]{sambusetti1999minimal}
Andrea Sambusetti.
\newblock Minimal entropy and simplicial volume.
\newblock {\em manuscripta mathematica}, 99(4):541--560, 1999.

\bibitem[SS09]{suarez2009minimal}
P.~Su{\'a}rez-Serrato.
\newblock Minimal entropy and geometric decompositions in dimension four.
\newblock {\em Algebr.\ Geom.\ Topol.}, 9(1):365--395, 2009.

\bibitem[Ste20]{Steinhuber}
M.~Steinhuber.
\newblock Integral foliated simplicial volume of mapping tori.
\newblock Master thesis, Universit{\"a}t Regensburg, 2020.

\bibitem[{\v{S}}va55]{vsvarc1955volume}
A.~S. {\v{S}}varc.
\newblock A volume invariant of coverings.
\newblock In {\em Dokl. Akad. Nauk SSSR (NS)}, volume 105, pages 32--34, 1955.

\end{thebibliography}

\end{document}